\documentclass[12pt]{article}
\usepackage[
    letterpaper,
    top=2cm,
    bottom=2cm,
    left=3cm,
    right=3cm,
    marginparwidth=1.75cm
]{geometry}

\usepackage[giveninits=true, backend=biber, maxbibnames=99]{biblatex}
\usepackage[english]{babel}

\usepackage[utf8]{inputenc}
\usepackage[T1]{fontenc}

\usepackage{amsmath}

\usepackage{amssymb}
\usepackage{amsthm}
\usepackage{booktabs}
\usepackage{indentfirst}

\newtheorem{proposition}{Proposition}
\usepackage{graphicx}
\usepackage{caption}
\usepackage{subfigure}
\usepackage{subcaption}
\usepackage{algorithm}
\usepackage{algpseudocode}
\usepackage{comment}
\usepackage[colorlinks=true, allcolors=blue]{hyperref}

\newcommand{\br}{\begin{eqnarray}}
\newcommand{\er}{\end{eqnarray}}

\title{A Stochastic Interacting Particle-Field Algorithm for a Haptotaxis Advection-Diffusion System  Modeling Cancer Cell Invasion }

\author{Boyi Hu\thanks{Department of Mathematics, The University of Hong Kong, Pokfulam Road, Hong Kong SAR, P.R.China. huby22@connect.hku.hk.}
\and 
Zhongjian Wang \thanks{Division of Mathematical Sciences, School of Physical and Mathematical Sciences, Nanyang Technological University, 21 Nanyang Link, Singapore 637371. zhongjian.wang@ntu.edu.sg}
\and
Jack Xin \thanks{Department of Mathematics, University of California at Irvine, Irvine, CA 92697, USA. jack.xin@uci.edu}
\and
Zhiwen Zhang\thanks{Corresponding author. Department of Mathematics, The University of Hong Kong, Pokfulam Road, Hong Kong SAR, P.R.China. Materials Innovation Institute for Life Sciences and Energy (MILES), HKU-SIRI, Shenzhen, P.R. China. zhangzw@hku.hk.}
 }

\begin{document}
\maketitle

\begin{abstract}
The investigation of tumor invasion and metastasis dynamics is crucial for advancements in cancer biology and treatment. Many mathematical models have been developed to study the invasion of host tissue by tumor cells. In this paper, we develop a novel stochastic interacting particle-field (SIPF) algorithm that accurately simulates the cancer cell invasion process within the haptotaxis advection-diffusion (HAD) system. Our approach approximates solutions using empirical measures of particle interactions, combined with a smoother field variable - the extracellular matrix concentration (ECM) - computed by the spectral method. We derive a one-step time recursion for both the positions of stochastic particles and the field variable using the implicit Euler discretization, which is based on the explicit Green's function of an elliptic operator characterized by the Laplacian minus a positive constant. Our numerical experiments demonstrate the superior performance of the proposed algorithm, especially in computing cancer cell growth with thin free boundaries in three-dimensional (3D) space. Numerical results show that the SIPF algorithm is mesh-free, self-adaptive, and low-cost. Moreover, it is more accurate and efficient than traditional numerical techniques such as the finite difference method (FDM) and spectral methods. \\

\end{abstract}

\section{Introduction}

The prevalence of cancer globally has increased significantly, making it the second leading cause of death following cardiovascular diseases \cite{siegel2024cancer}. The dynamics of tumor invasion and metastasis are crucial research areas within cancer biology and treatment. Since the 1970s, various mathematical models have been developed to analyze the different phases of solid tumor growth, both in temporal and spatio-temporal contexts \cite{chaplain1996avascular}. A significant amount of empirical data on the growth dynamics of avascular tumors has been integrated into mathematical models that utilize various growth laws, including Gompertzian, logistic, and exponential growth \cite{retsky1990gompertzian}. Additionally, stochastic growth models have been utilized to simulate the invasion of tumor cells, providing insights into the functional implications of histological patterns \cite{smolle1993computer}.

Mathematical modeling is a powerful tool for unraveling the complexities of biological processes, providing insights that inform both experimental and clinical strategies. The field of cancer modeling has benefited from various approaches, from mechanistic models that explore the detailed mechanisms of diseases to data-driven models that facilitate clinical decision-making \cite{bekisz2020cancer}. Specifically, in the area of tumor-induced angiogenesis, researchers have developed both continuous and discrete mathematical models to simulate the formation of capillary networks triggered by tumor angiogenic factors. These models effectively integrate critical interactions between endothelial cells and the ECM \cite{ANDERSON1998857}. Cancer cell invasion of tissue is a complex process where cell migration through the ECM, facilitated by the secretion of degradative enzymes, plays a central role \cite{RAMISCONDE2008533}. This invasion is modeled using a system of partial differential equations (PDEs) that capture the dynamics involving tumor cells, the ECM, and matrix-degrading enzymes (MDEs), highlighting the intricate biological interactions essential for tissue invasion \cite{anderson2000mathematical}.

Additionally, fractional mathematical models have been introduced to better understand the intricate dynamics among tumor cells, matrix degradation, and enzyme production, employing sophisticated analytical techniques such as the $q$-homotopy analysis transform method \cite{veeresha2021regarding}. Models employing stochastic differential equations (SDE) have been formulated to capture the stochastic behaviors of cancer cell migration and invasion, specifically addressing the variability in diffusion processes within the context of PDE  \cite{katsaounis2023stochastic}. Furthermore, the global behavior of solutions to models of tumor invasion, which emphasize the critical role of ECM concentration, has undergone thorough analysis, providing a detailed understanding of the invasion process and its interactions with various biological factors \cite{tao2007global, marciniak2010boundedness, lictcanu2010asymptotic}. The integration of these models forms a comprehensive framework of mathematical and computational approaches that significantly enhance our understanding of cancer cell invasion and metastasis. This collective body of work lays a solid foundation for the innovation of therapeutic strategies aimed at combating cancer effectively. 

Given the substantial theoretical progress made using the Lagrangian perspective in models analogous to cancer invasion, it is compelling to consider applying this framework to numerically solve problems related to cancer invasion. A convergent particle method was derived for fully parabolic chemotaxis equations \cite{stevens2000derivation}. The study \cite{stevens2000stochastic} utilized cellular automaton simulations based on a reinforced self-attracting random walk for a single particle in 1D. Building on this framework, \cite{stevens1997aggregation} expanded the scope to derive more general chemotaxis systems from similar reinforced random walks, and further analyzed the qualitative behavior of these systems, providing deeper insights into the dynamics and implications of such models in higher dimensions and more complex scenarios. A random particle blob method is shown to converge for the parabolic-elliptic Keller-Segel (KS) system when the macroscopic mean-field equation allows for a global weak solution \cite{liu2019propagation,liu2017random}. The success of this method greatly hinges on an in-depth understanding of the nonlinear mean-field equation, rather than the complexities of the multi-particle Markov process involved \cite{mischler2013kac}. A deep-learning study of chemotaxis and aggregation in 3D laminar and chaotic flows is initiated in \cite{wang2024deepparticle} with a kernel regularization technique for particle dynamics.

In this paper, we introduce a SIPF algorithm to compute the cancer cell invasion process, as proposed \cite{wang2023novel} for the fully parabolic KS system. Our method takes into account the coupled stochastic particle and field evolution, where the corresponding field represents the MDE concentration within the system. This approach enables self-adaptive simulations that effectively handle potential singularities or free boundaries. In our SIPF algorithm, we model the density of active particles using empirical particle representation, which entails a summation of delta functions centered at the particle positions. Furthermore, we discretize the MDE concentration field using the spectral method instead of FDM, as suggested by \cite{fatkullin2012study}. This choice is motivated by the fact that the field tends to be smoother than the density. Specifically, the ECM density is updated through an explicit Euler scheme applied to its Fourier coefficients, leveraging the convolution theorem. 

To validate the efficacy of our method, we carried out numerical experiments in three dimensions (3D). It is worth mentioning that pseudo-spectral methods have been effectively employed in computing nearly singular solutions of the 3D Euler equations \cite{hou2007computing}. Additionally, the adaptive moving mesh method has been developed to investigate finite-time blowup in the 3D axisymmetric Euler equations \cite{luo2014toward}. These approaches are high-resolution methods for resolving nearly singular phenomena in the 3D Euler equations. On the other hand, 
it is important to acknowledge that the implementation of pseudo-spectral methods for 3D problems requires significant computational resources, and the adaptive moving mesh method requires intricate design and advanced coding capabilities. 
  

In contrast, SIPF is a simple-to-program and efficient low-cost method for 3D computations. We show the effectiveness of the SIPF algorithm for simulating the HAD system modeling cancer cell invasion. The algorithm operates recursively without relying on historical data and computes the field variable (concentration) using Fast Fourier Transforms (FFT). The concentration field variable is smoother than the particle density, allowing FFT to work efficiently with only a few dozen Fourier modes. The spreading phenomenon of the particle density is accurately captured using 10,000 particles. Traditional implicit FDM is not only time-consuming but also suffers from poor precision in 3D numerical simulations, which results in an inaccurate representation of tumor invasion as it occurs in reality. 

The rest of the paper is organized as follows. In Section 2, we provide a concise overview of the cancer cell invasion HAD system, with a discussion of related theoretical analyses. We introduce our SIPF algorithms, which streamline a theoretically equivalent yet computationally intensive method (involving history-dependent parabolic kernel functions) into practical recursive computations. In Section 3, we present numerical results to illustrate the efficacy of SIPF within 3D cancer cell invasion models. It outperforms FDM in terms of both computational runtime and accuracy, especially when the diffusion coefficient of cancer cell density is small. The paper concludes with Section 4, where we summarize our results and discuss potential avenues for future research.

\section{Cancer Cell Invasion HAD System}
 
Consider the following Cancer Cell Invasion HAD System (see Eq.(5) in \cite{anderson2000mathematical}), which describes the interactions among tumor cells, ECM, and MDEs:

\begin{equation}\label{cancer system}
\begin{aligned}
\rho_t &= d_n \Delta \rho - \gamma \nabla \cdot (\rho \nabla f), \\
f_t &= -\eta m f, \\
m_t &= d_m \Delta m - \beta m + \alpha \rho.
\end{aligned}
\end{equation}

The system is defined with physical and biological parameters \((d_n, \gamma, \eta, d_m, \alpha, \beta) > 0\) on a compact subset \(\Omega\) of \(\mathbb{R}^d\) (where \(d = 2, 3\)), and zero flux boundary conditions (refer to Eq.(6)-(7) in  \cite{anderson2000mathematical}, with \(\partial_n\) representing the outward normal derivative):
\begin{equation}\label{boundary condition}
d_n \partial_n\rho - \gamma \rho \partial_n f = 0, \quad \partial_nm = 0; \quad \text{on } \partial\Omega.
\end{equation}
Given that $f = f_0 \exp\left(-\eta \int_0^t m(x; \tau) \, d\tau\right)$, 
we have $\partial_n f = 0$ on $\Omega$ if $(\partial_n f_0, \partial_n m) = 0$ on $\partial \Omega$. Hence, \eqref{boundary condition} is replaced by the following zero Neumann boundary conditions:
\begin{equation}\label{zNbc boundary condition}
\partial_n(\rho, m) = 0, \quad \forall t \geq 0 \text{ on } \partial\Omega; \quad \partial_n f_0= 0 \text{ on } \partial\Omega.
\end{equation}

The variables $\rho$, $m$, and $f$ in \eqref{cancer system} are functions of both the spatial variable $x$ and time $t$. The first equation in \eqref{cancer system} governs tumor cell motion, with $\rho$ representing the tumor cell density. Our model specifically focuses on the interactions between cells and the ECM, examining their impact on tumor cell migration without incorporating cell proliferation. We choose $d_n$ to be a constant, representing the tumor cell random motility coefficient, rather than a function of either the MDE or ECM concentration. 

The second equation in \eqref{cancer system} models this degradation process, with $f$ representing the ECM density and $\delta$ being a positive constant. We posit that the MDEs degrade the ECM when they come into contact. The tumor cells generate (or activate) the MDEs, which diffuse throughout the surrounding tissue. These active MDEs subsequently undergo some form of decay (either passive or active).

The third equation in \eqref{cancer system} governs the evolution of MDE concentration, with $m$ representing the MDE concentration and the positive constant $d_m$ for the MDE diffusion coefficient. For the sake of simplicity, the model stipulates a straightforward, proportional connection between the density of tumor cells and the amount of active MDEs present in the surrounding tissues. This linear association is presumed to hold independent of the quantity of enzyme precursors secreted or the existence of endogenous inhibitors. Consequently, we initially take $\alpha \rho$ to represent the MDE production by the tumor cells and $\beta m$ to signify natural decay, respectively.

\subsection{Global Well-Possedness of Nonnegative Solutions}

The paper \cite{lictcanu2010asymptotic} provides some theoretical analyses for the general reaction-advection-diffusion model, which has also been discussed in \cite{marciniak2010boundedness} and \cite{gopika2023residual}.
We follow the line of proof in \cite{lictcanu2010asymptotic}. By setting $\alpha_2=0$, $g(v)\equiv 1$, and $\chi(v)\equiv 1$ in \cite{lictcanu2010asymptotic}, the system in \cite{lictcanu2010asymptotic} simplifies to the same as the system \eqref{cancer system}. First, we nondimensionalize the system \eqref{cancer system}, and it becomes:

\begin{equation}\label{nondimensionalizing cancer system}
\begin{aligned}
\rho_t &= \Delta \rho - \nabla \cdot (\rho \nabla f), \\
f_t &= -m f, \\
m_t &= d\Delta m - \lambda m + \rho, \\
\partial_n \rho - \rho \partial_n f &= \partial_n m = 0, & x \in \partial\Omega,
\end{aligned}
\end{equation}
where  
$d = \frac{d_m}{d_n}$ and $\lambda = \frac{\beta}{\alpha}$. 

The positive self-adjoint operator $A$ is defined as:
\[
A := -a\Delta + b
\]
where the domain $D(A)$ is given by:
\[
D(A) := \{ \rho \in W^{2,p}(\Omega) : \frac{\partial \rho}{\partial n} = 0 \text{ on } \partial\Omega \}
\]
Here, $p \in (1,+\infty)$, $a>0,b>0$. 

For \(0 \leq \theta \leq 1\), the fractional powers of the operator $A$, denoted as $A^\theta$, map the space $X_p^\theta$ to $L^p(\Omega)$. The space $X_p^\theta$ is equipped with the graph norm:
\[
\|u\|_{X_p^\theta} = \|A^\theta u\|_{L^p(\Omega)}.
\]

\begin{proposition}[Theorem 3.3, \cite{lictcanu2010asymptotic}] 
Let $\Omega \subset \mathbb{R}^N$, $N \geq 1$ be a domain with $C^2$ boundary and $p > N$. Given the non-negative initial value $(\rho_0, f_0, m_0) \in W^{1,p}(\Omega) \times W^{1,\infty}(\Omega) \times X_{p}^{\theta}$, $\theta \in \left(\frac{N+p}{2p}, 1\right)$, there exists $T > 0$ (depending only on $\|\rho_0\|_{W^{1,p}(\Omega)}$, $\|f_0\|_{W^{1,\infty}(\Omega)}$ and $\|m_0\|_{X_{p}^{\theta}}$) such that the system \eqref{nondimensionalizing cancer system} has a unique non-negative solution $(\rho, f, m)$ defined on an interval $[0, T) \subset \mathbb{R}$ and
\begin{align*}
\rho &\in C\left([0, T); W^{1,q}(\Omega)\right) \cap C\left((0, T); W^{1,\infty}(\Omega)\right) \cap C^1\left((0, T); W^{1,q}(\Omega)\right), \\
f &\in C\left([0, T); W^{1,\infty}(\Omega)\right) \cap C^1\left((0, T); W^{1,\infty}(\Omega)\right), \\
m &\in C\left([0, T); X_{p}^{\theta}\right) \cap C\left((0, T); W^{2,p}(\Omega)\right) \cap C^1\left((0, T); X_{p}^{\theta}\right).
\end{align*}
Moreover, the solution depends continuously on the initial data.
\end{proposition}

Proposition 1 implies the local existence, uniqueness, and non-negativity of the solution to the non-dimensionalized cancer system \eqref{nondimensionalizing cancer system} when the initial values $(\rho_0, f_0, m_0)$ are non-negative. Based on the lemmas and propositions in Section 4 of \cite{lictcanu2010asymptotic}, we can derive the following proposition, which demonstrates the global existence of the solution in $\Omega \subset \mathbb{R}^3$.

\begin{proposition}
Let \( \Omega \subset \mathbb{R}^3 \) be a domain with a smooth boundary. Given the non-negative initial value \( (\rho_0, f_0, m_0) \in L^{\infty}(\Omega) \times W^{1,\infty}(\Omega) \times X_p^\theta(\Omega) \), where \( p > 3 \), \( \theta \in \left( \frac{3+p}{2p}, 1 \right) \), for all \( t \in [0, T) \), there exists a constant \(C_p\) independent on time such that the solution $(\rho, f, m)$ to the system \eqref{nondimensionalizing cancer system} satisfies that \[\|\rho(\cdot, t)\|_{L^{\infty}(\Omega)} + \|f(\cdot, t)\|_{L^p(\Omega)} + \|m(\cdot, t)\|_{X_p^\theta(\Omega)} \leq C_p.\] Since the solution \((\rho, f, m)\) is uniformly bounded in \(L^\infty(\Omega)\) for all \(t \in [0, T)\), the local solution can be extended to a global solution.
\end{proposition} 

\subsection{Integral Identities}
By integrating the equation of MDE concentration (second equation in \eqref{cancer system}) over the spatial domain, we have: 
\begin{equation}
\int_{\Omega} \Delta m \, dx= \int_{\partial \Omega} \nabla m \cdot \mathbf{n} \, dS = 0, \;\; \text{due}
\; \text{to}\;\; 
\nabla m \cdot \mathbf{n} = 0 \quad \text{on} \quad \partial \Omega,
\end{equation}
\begin{equation}\label{integral m}
\frac{d}{dt} \int_{\Omega} m \, dx = -\beta \int_{\Omega} m\, dx  + \alpha \int_{\Omega} \rho\, dx,
\end{equation}
where $\int_{\Omega} \rho\, dx$ is conserved in time given the conservation form of the $\rho$ equation (\ref{nondimensionalizing cancer system}). Hence $\int_{\Omega}\, m\, dx$ can be integrated in closed analytic form given its initial value. 
In terms of $m$, $f$ satisfies: 
\begin{equation}\label{f0 and m}
f = f_0 \exp\left(-\eta \int_{0}^{t} m(x; \tau) d\tau\right).
\end{equation}
Here $f_0$ refers to the initial value of $f$ at $t=0$. 
Taking the logarithm and integrating the equation \eqref{f0 and m} over space, we get 
\begin{equation}
\int_{\Omega}  \ln f\, dx = \int_{\Omega}  \ln f_0\, dx-\eta \int_{\Omega}\int_{0}^{t} m(x; \tau) d\tau dx.
\end{equation}
From the well-posedness of the system (\ref{nondimensionalizing cancer system}), \( m(x; \tau) \) is integrable in 
\( x \) and \( \tau \) and the integral of \( m(x; \tau) \) over \( \Omega \times [0, t] \) is finite, according to the Fubini's theorem, the order of integration can be interchanged, which implies that
\begin{equation}\label{lnf}
\int_{\Omega}  \ln f\, dx = \int_{\Omega}  \ln f_0\, dx-\eta \int_{0}^{t}\int_{\Omega} m(x; \tau)\, dx \,d\tau.
\end{equation}
So $\int_{\Omega}  \ln f\, dx$ is also in closed analytical form via time integration of 
$\int_\Omega \, m \, dx$. These identities will be used to validate numerical approximations later in \eqref{m relative error}.

\subsection{SIPF Algorithm for Cancer Cell Invasion HAD System}
To solve the model of invasion of host tissue by tumor cells, we adopt a similar SIPF algorithm  \cite{wang2023novel} for KS systems. Let us partition $[0, T]$ by temporal grid points
$\left\{t_n\right\}_{n=0:n_{T}}$ with $t_0 = 0$ and $t_{n_{T}} = T$, and approximate the density $\rho$ by particles as follows:
\begin{equation}\label{particle rho}
 \rho_t \approx \frac{M_0}{P}\sum_{j=1}^{P}\delta(\mathbf{x}-X_t^{P}), \;\;\; P \gg 1,
\end{equation}
where $M_0$ is the conserved total mass (integral of $\rho$), and $P$ is the number of particles. 
 
We restrict the domain $\Omega$ to be $[0, L\pi]^d$ due to our application of the Fourier transform. For the MDE concentration $m$, the method to treat the chemical concentration field $c$ in \cite{wang2023novel} is used. We approximate $m(\mathbf{x}, t)$ by Fourier series:
\begin{equation}\label{m Fourier}
 \sum_{j,m,l \in \mathcal{H}}\alpha_{t;j,m,l}\exp(i2 \pi jx_1/L)\exp(i2 \pi mx_1/L)\exp(i2 \pi lx_1/L)
\end{equation}
where $\mathcal{H}$ denotes the index set
\begin{equation}
 \{(j,m,l) \in \mathbb{N}^3: |j|, |m|, |l| \leq \frac{H}{2}\},
\end{equation}
and $i=\sqrt{-1}$.

Then at \( t_0 = 0 \), we generate \( P \) empirical samples \( \{X_{0}^{p}\}_{p=1}^{P} \) according to the initial condition of \( \rho_0 \) and set up \( \alpha_{0;j,m,l} \) by the Fourier series of \( m_0 \).

For ease of presenting our algorithm, with a slight abuse of notation, we use \( \rho_n = \frac{M_0}{P} \sum_{p=1}^{P} \delta(\mathbf{x} - X_{n}^{p}) \), \[ m_n = \sum_{j,m,l \in \mathcal{H}} \alpha_{n;j,m,l} \exp(i2\pi j x_1/L) \exp(i2\pi m x_2/L) \exp(i2\pi l x_3/L) \] and \[f_n =  \sum_{j,m,l \in \mathcal{H}} \beta_{n;j,m,l} \exp(i2\pi j x_1/L) \exp(i2\pi m x_2/L) \exp(i2\pi l x_3/L) \]
to represent tumor cell density \( \rho \), MDE concentration \( m \) and ECM density \( f \) at time \( t_n \).
Considering time stepping system \eqref{cancer system} from $t_n$ to $t_{n+1}$, with $\rho_n$, $f_{n-1}$ and $m_{n-1}$ known, our algorithm, inspired by the operator splitting technique, consists of three sub-steps: updating MDE concentration $m$, updating ECM density $f$ and updating tumor cell density $\rho$.

\vspace{12pt}
\noindent
\textit{Updating MDE concentration \( m \)}. Let \( \delta t = t_{n+1} - t_n > 0 \) be the time step. We discretize the \( m \) equation of \eqref{cancer system} in time by an implicit Euler scheme:
\begin{equation}\label{implicit m}
 \frac{m_n-m_{n-1}}{\delta t}=d_m \cdot \Delta m_n - \beta m_n + \alpha \rho_n.
\end{equation}
From \eqref{implicit m}, we obtain the explicit formula for $m_n$ as:
\begin{equation}
 (\Delta-\frac{\beta}{d_m}-\frac{1}{d_m \cdot \delta t})\cdot m_n=-(\frac{1}{d_m \cdot \delta t}\cdot m_{n-1}+\frac{\alpha}{d_m}\cdot \rho_n). 
\end{equation}
It follows that:
\begin{equation}\label{convolution m}
 m_n=m(\mathbf{x},t_{n})=-\mathcal{K}_{\delta t}\ast(\frac{m_{n-1}}{d_m \cdot \delta t}+\frac{\alpha}{d_m}\cdot\rho_{n})=-\mathcal{K}_{\delta t}\ast(\frac{m(\mathbf{x},t_{n-1})}{d_m \cdot \delta t}+\frac{\alpha}{d_m}\cdot\rho(\mathbf{x},t_{n})).
\end{equation}
where $\ast$ is spatial convolution operator, and $\mathcal{K}_{\delta t}$ is the Green’s function of the operator $(\Delta-\frac{1}{d_m \cdot \delta t})$. 
\br
\nabla_{\mathbf{x}}m_n&=&\nabla_{\mathbf{x}}m(\mathbf{x},t_{n})=-\nabla_{\mathbf{x}}\mathcal{K}_{\delta t}\ast(\frac{m_{n-1}}{d_m \cdot \delta t}+\frac{\alpha}{d_m}\cdot\rho_{n})\nonumber \\
 &=& -\nabla_{\mathbf{x}}\mathcal{K}_{\delta t}\ast(\frac{m(\mathbf{x},t_{n-1})}{d_m \cdot \delta t}+\frac{\alpha}{d_m}\cdot\rho(\mathbf{x},t_{n})).
\er
In $\mathbb{R}^3$, the Green’s function $\mathcal{K}_{\delta t}$ reads as follows:
\begin{equation}
 \mathcal{K}_{\delta t}=\mathcal{K}_{\delta t}(\mathbf{x})=-\frac{e^{-|\zeta|\mathbf{x}}}{4\pi |\mathbf{x}|},\;\; \zeta^2=\frac{\beta}{d_m}+\frac{1}{d_m \cdot \delta t}.
\end{equation}
The Green’s function admits a closed-form Fourier transform,
\begin{equation}\label{fourier transform m}
 \mathcal{F}\mathcal{K}_{\delta t}(\omega)=-\frac{1}{|\omega|^2+\zeta^2}.
\end{equation}
For the term \( -\mathcal{K}_{\delta t} * m_{n-1} \) in \eqref{convolution m}, by Eq.\eqref{fourier transform m} it is equivalent to modify Fourier coefficients \( \alpha_{n;j,m,l} \) to
\[ \frac{\alpha_{n;j,m,l}}{4\pi^2 j^2 / L^2 + 4\pi^2 m^2 / L^2 + 4\pi^2 l^2 / L^2 + \zeta^2}. \]
For the second term \( \mathcal{K}_{\delta t} * \rho \), we first approximate \( \mathcal{K}_{\delta t} \) with a cosine series expansion, then according to the particle representation of \( \rho \) in \eqref{particle rho},
\[ (\mathcal{K}_{\delta t} * \rho)_{j,m,l} \approx \frac{M_0}{P} \sum_{p=1}^{P} \frac{\exp\left(-2\pi j X_{n,1}^p/L - 2\pi m X_{n,2}^p/L - 2\pi l X_{n,3}^p/L\right)(-1)^{j+m+l}} {\left(4\pi^2 j^2/L^2 + 4\pi^2 m^2/L^2 + 4\pi^2 l^2/L^2 + \zeta^2\right)}. \]
Finally, we summarize the one-step update of Fourier coefficients of MDE concentration \( m \) in Alg.\ref{update MDE}.
\begin{algorithm}
\caption{One step update of MDE concentration in SIPF}
\label{update MDE}
\begin{algorithmic}[1]
\Require Distribution $\rho_n$ represented by empirical samples $X_n$, initial MDE concentration $m_{n-1}$ represented by Fourier coefficients $\alpha_{n-1}$
\ForAll{$(j, m, l) \in \mathcal{H}$}
    \State $\alpha_{n;j,m,l} \gets \frac{\alpha_{n-1;j,m,l}}{d_m \cdot \delta t(4\pi^2j^2/L^2 + 4\pi^2m^2/L^2 + 4\pi^2l^2/L^2 + \zeta^2)}$
    \State $F_{j,m,l} \gets 0$
    \For{$p = 1$ to $P$}
        \State $F_{j,m,l} \gets F_{j,m,l} + \exp\left(-2\pi jX_{n;1}^{p}/L - 2\pi mX_{n;2}^{p}/L - 2\pi lX_{n;3}^{p}/L\right)$
    \EndFor
    \State $F_{j,m,l} \gets F_{j,m,l} \frac{(-1)^{j+m+l}}{4\pi^2j^2/L^2 + 4\pi^2m^2/L^2 + 4\pi^2l^2/L^2 + \beta^2}\ast\frac{M}{P}$
\EndFor
\State $\alpha_n \gets \alpha_n - \frac{\alpha}{d_m}\cdot F$
\Ensure Updated MDE concentration field from input $m_{n-1}$ to $m_n$ via $\alpha_n$.
\end{algorithmic}
\end{algorithm}

\vspace{12pt}
\noindent
\textit{Updating ECM density \( f \)}.
$f(\mathbf{x}, t)$ has an series representation:
\begin{equation}
   \sum_{j,m,l \in \mathcal{H}}\beta_{t;j,m,l}\exp(i2 \pi jx_1/L)\exp(i2 \pi mx_2/L)\exp(i2 \pi lx_3/L).
   \end{equation}
We discretize the \( f \) equation of \eqref{cancer system} in time by an explicit Euler scheme:
\begin{equation}
f(\mathbf{x},t_{n+1})=f(\mathbf{x},t_n)-\eta \cdot m(\mathbf{x},t_n)f(\mathbf{x},t_n)\delta t
\end{equation}
For $f(\mathbf{x},t_{n+1})$, according to the convolution theorem, it is equivalent to modify Fourier coefficients \( \beta_{n;j,m,l} \) to\[\beta_{n;j,m,l} - \eta\sum_{j',m',l' \in \mathcal{H}} \alpha_{n;j',m',l'} \beta_{n;j-j',m-m',l-l'} \, \delta t. \]
It follows that:
\begin{equation}\label{grad f}
\nabla_{\mathbf{x}} f(\mathbf{x},t_{n+1})=\nabla_{\mathbf{x}} f(\mathbf{x},t_n)-\eta f(\mathbf{x},t_n)\delta t\nabla_{\mathbf{x}} m(\mathbf{x},t_n)-\eta m(\mathbf{x},t_n)\delta t\nabla_{\mathbf{x}} f(\mathbf{x},t_n).
\end{equation}
Finally, we summarize the one-step update of Fourier coefficients of ECM density \( f \) in Alg.\ref{update ECM}.

\begin{algorithm}[H] 
\caption{One step update of ECM density in SIPF}
\label{update ECM}
\begin{algorithmic}[1] 
\Require Initial MDE concentration \( m_{n-1} \) represented by Fourier coefficients \( \alpha_{n-1} \), initial ECM density \( f_{n-1} \) represented by Fourier coefficients \( \beta_{n-1} \)
\ForAll{\( (j, m, l) \in \mathcal{H} \)}
   \State \( F_{j,m,l} \gets 0 \)
\ForAll{\( (j', m', l') \in \mathcal{H} \)}
   \State \( F_{j,m,l} \gets F_{j,m,l} - \eta \alpha_{n-1;j',m',l'}\beta_{n-1;j-j',m-m',l-l'}\delta t \)
\EndFor
   \State \( \beta_{n;j,m,l} \gets \beta_{n-1;j,m,l} - F_{j,m,l} \)
\EndFor
\Ensure Updated ECM density from input \( f_{n-1} \) to \( f_n \) via \( \beta_n \).
\end{algorithmic}
\end{algorithm}

\vspace{12pt}
\noindent
\textit{Updating density of active particles \( \rho \)}. 
After updating the ECM density $f$, we can update the density $\rho$. The empirical particle system converging to density $\rho$ reads:
\begin{equation}\label{rho SDE}
    dX_p = \gamma \nabla_\mathbf{x} f(X_{t}^p, t) \, dt + \sqrt{2d_n} \, dW_p, \quad p = 1, \ldots, P,
\end{equation}
where $W_p$ are independent standard Brownian motions in $\mathbb{R}^d$.

In the one-step update of density $\rho_n$ represented
by particles $\left\{X_n^{p}\right\}_{p=1:P}$, we apply Euler-Maruyama scheme to solve the SDE \eqref{rho SDE}:

\begin{equation}\label{X SDE}
 X_{n+1}^{p}=X_{n}^{p}+\gamma \nabla_{\mathbf{x}}f(X_n^{p},t_n)\cdot \delta t+\sqrt{2d_n\delta t}N_n^{p}
\end{equation}
where ${N_n^{p}}'s$ are i.i.d. standard normal distributions with respect to the Brownian paths
in the SDE formulation \eqref{rho SDE}. 
For $n$ > 1, substituting \eqref{grad f} in \eqref{X SDE} gives:
\begin{equation}
   \begin{split}
     X_{n+1}^{p} &= X_{n}^{p} + \gamma \big( \nabla_{\mathbf{x}}f(\mathbf{x},t_{n-1}) - \eta \nabla_{\mathbf{x}}f(\mathbf{x},t_{n-1})\cdot m(\mathbf{x},t_{n-1})\delta t \\
    &- \eta \nabla_{\mathbf{x}}m(\mathbf{x},t_{n-1})\cdot f(\mathbf{x},t_{n-1})\delta t \big) \Big|_{\mathbf{x} = X_{n}^{p}}\delta t + \sqrt{2d_n\delta t}N_n^{p}.
   \end{split}
\end{equation}

In such particle formulation, the computation of spacial convolution is slightly different from the one in the update of $m$, namely \eqref{convolution m}.
\begin{equation}
   \nabla_{\mathbf{x}}m(\mathbf{x},t_{n-1})=-\nabla_{\mathbf{x}}\mathcal{K}_{\delta t}\ast(\frac{m(\mathbf{x},t_{n-2})}{d_m \cdot \delta t}+\frac{\alpha}{d_m}\cdot \rho(\mathbf{x},t_{n-1})).
\end{equation}

For $\nabla_{\mathbf{x}} \mathcal{K}_{\delta t} \ast m_{n-2}(X_n^p)$, to avoid the singular points of $\nabla_{\mathbf{x}} \mathcal{K}_{\delta t}$, we evaluate the integral with the quadrature points that are away from 0. 
Precisely, denote the standard quadrature point in $\Omega$ with
\begin{equation}
   x_{j,m,l} = \left(\frac{jL}{H}, \frac{mL}{H}, \frac{jL}{H}\right),
\end{equation} 
where $j, m, l$ are integers ranging from $-\frac{H}{2}$ to $\frac{H}{2} - 1$. 

We evaluate $\nabla_\mathbf{x} \mathcal{K}_{\delta t}$ at $\{X_n^p + \overline{X}_n^p - x_{j,m,l}\}_{j,m,l}$ where a small spatial shift $\overline{X}_n^p$ is defined as ${\overline{X}_n^p = \frac{H}{2L} + \left\lfloor \frac{X_n^p}{H/L} \right\rfloor \frac{H}{L} - X_n^p}$, and $m$ at $\{x_{j,m,l} - \overline{X}_n^p\}_{j,m,l}$. The latter one is computed by inverse Fourier transform of shifted coefficients, with $\alpha_{j,m,l}$ modified to 
\[ \alpha_{j,m,l} \exp\left(-i2\pi j \overline{X}_{n;1}^p/L- i2\pi m \overline{X}_{n;2}^p/L - i2\pi l \overline{X}_{n;3}^p/L\right)\] 
where $\overline{X}_{n;i}^p$ denotes the $i$-th component of $\overline{X}_n^p$. 

The term $\nabla_\mathbf{x} K_{\delta t} \ast \rho(X_n^p, t_{n-1})$ is straightforward thanks to the particle representation of $\rho(X_n^p, t_{n-1})$ in \eqref{particle rho}:
\begin{equation}
\nabla_\mathbf{x} \mathcal{K}_{\delta t} \ast \rho_{n-1}(X_n^p) = \int \mathcal{K}_{\delta t}(X_n^p - y)\rho(y) \, dy \approx \sum_{\substack{q=1}}^P \frac{M}{P} \mathcal{K}_{\delta t}(X_n^p - X_{n-1}^q).
\end{equation}

Finally, we summarize the one-step update of Fourier coefficients of tumor density \( \rho \) in Alg.\ref{update tumor cell density}.

\begin{algorithm}[H]
\caption{One step update of tumor cell density in SIPF}
\label{update tumor cell density}
\begin{algorithmic}[1]
\State \textbf{Data}: Distribution $\rho_n$ represented by empirical samples $X_n$, input MDE concentration $m_{n-1}$ represented by Fourier coefficients $\alpha_{n-1}$, ECM density $f_{n-1}$ represented by Fourier coefficients $\beta_{n-1}$, $\rho_{n-1}$ represented by empirical samples $X_{n-1}$,  MDE concentration $m_{n-2}$ represented by Fourier coefficients $\alpha_{n-2}$
\For{$p = 1$ to $P$}
\State $X_{n+1}^p \gets X_n^p + \sqrt{2 d_n\delta t}N$ \Comment{where $N$ is a standard normal distribution}
\State $\bar{X}_n^p \gets \frac{H}{2L} + \left\lfloor \frac{X_n^p}{H/L} \right\rfloor \frac{H}{L} - X_n^p$
\State $\nabla_{\mathbf{x}}f(X_n^p,t_{n-1}) \gets 0$; $f(X_n^p,t_{n-1}) \gets 0$; $m(X_n^p,t_{n-1}) \gets 0$
\ForAll{$(j, m, l) \in \mathcal{H}$}
\State $F_{j,m,l} \gets \nabla_\mathbf{x} \mathcal{K}_{\epsilon,\delta t}(X_n^p + \bar{X}_n^p - x_{j,m,l})$ \Comment{$x_{j,m,l} = \left(jL/H, mL/H, lL/H\right)$}
\State $G_{j,m,l} \gets \alpha_{n-2;j,m,l} \exp\left(-2\pi j\bar{X}_{n;1}^{p}/L - 2\pi m\bar{X}_{n;2}^{p}/L - 2\pi l\bar{X}_{n;3}^{p}/L\right)$
\State $\nabla_{\mathbf{x}}f(X_n^p,t_{n-1}) \gets \nabla_{\mathbf{x}}f(X_n^p,t_{n-1}) + \frac{i2 \pi}{L} \beta_{n-1;j,m,l}e^{i2\pi j X_{n;1}^p/L} e^{i2\pi m X_{n;2}^p/L} e^{i2\pi l X_{n;3}^p/L} \cdot (j, m, l)$ 
\State $f(X_n^p,t_{n-1}) \gets f(X_n^p,t_{n-1}) + \beta_{n-1;j,m,l}e^{i2\pi j X_{n;1}^p/L} e^{i2\pi m X_{n;2}^p/L} e^{i2\pi l X_{n;3}^p/L}$ 
\State $m(X_n^p,t_{n-1}) \gets m(X_n^p,t_{n-1}) + \alpha_{n-1;j,m,l}e^{i2\pi j X_{n;1}^p/L} e^{i2\pi m X_{n;2}^p/L} e^{i2\pi l X_{n;3}^p/L}$
\EndFor
\State $\hat{G} \gets \text{iFFT}(G)$
\For{$q = 1$ to $P$}
\State $\nabla_{\mathbf{x}}m(X_n^p,t_{n-1}) \gets  -\langle F, \hat{G}\rangle \frac{L^3}{H^3}/(d_m \cdot \delta t)- \frac{\alpha M}{d_m P} \mathcal{K}_{\delta t}(X_n^p - X_{n-1}^q)  $ \Comment{where $\langle\cdot, \cdot\rangle$ denote an inner product corresponding to $L^2(\Omega)$ quadrature}
\EndFor
\State $\nabla_{\mathbf{x}}f(X_n^p,t_{n}) \gets \nabla_{\mathbf{x}}f(X_n^p,t_{n-1})-\eta \nabla_{\mathbf{x}}f(\mathbf{x},t_{n-1})\cdot m(\mathbf{x},t_{n-1})\delta t-\eta \nabla_{\mathbf{x}}m(\mathbf{x},t_{n-1})\cdot f(\mathbf{x},t_{n-1})\delta t$ 
\State $X_{n+1}^p \gets X_{n+1}^p + \gamma \cdot \nabla_{\mathbf{x}}f(X_n^p,t_{n})\cdot \delta t$ 
\EndFor
\State \textbf{Result}: Output $\rho_{n+1}$ represented by updated $X_{n+1}$.
\end{algorithmic}
\end{algorithm}

\begin{algorithm}
\caption{Stochastic Interacting Particle-Field Method}
\label{the whole algorithm}
\begin{algorithmic}[1]
\State \textbf{Data}: Initial distribution $\rho_0$, initial MDE concentration $m_0$, initial ECM density $f_0$
\State Generate $P$ i.i.d samples following distribution $\rho_0$, $X_1, X_2, \ldots, X_P$.
\For{$p \gets 1$ to $P$}
    \State Compute $X_{1}^p$ by \eqref{X SDE}, with $f_0$
\EndFor
\State Compute $f_1$ by Alg.2 with $f_0$ and $m_0$.
\State Compute $m_1$ by Alg.1 with $m_0$ and $\rho_1 = \sum_{p=1}^{P} \frac{M}{P} \delta_{X_{1}^p}$.
\For{$\text{step } n \gets 2$ to $N = T/\delta t$}
    \State Compute $X_n$ by Alg.3 with $\rho_{n-1}$, $m_{n-1}$, $\rho_{n-2}$, $m_{n-2}$ and $f_{n-1}$
    \State Compute $f_n$ by Alg.2 with $f_{n-1}$ and $m_{n-1}$.
    \State Compute $m_n$ by Alg.1 with $m_{n-1}$ and $\rho_n =\sum_{p=1}^{P} \frac{M}{P} \delta_{X_{n}^p}$.
\EndFor
\end{algorithmic}
\end{algorithm}

\section{Numerical Experiments}
\begin{figure}[htbp]
	\centering
    \vspace{-0.2cm}
    \subfigtopskip=2pt
    \subfigure[t=0]{
	\begin{minipage}{0.49\linewidth}
		\centering
		\includegraphics[width=0.9\linewidth]{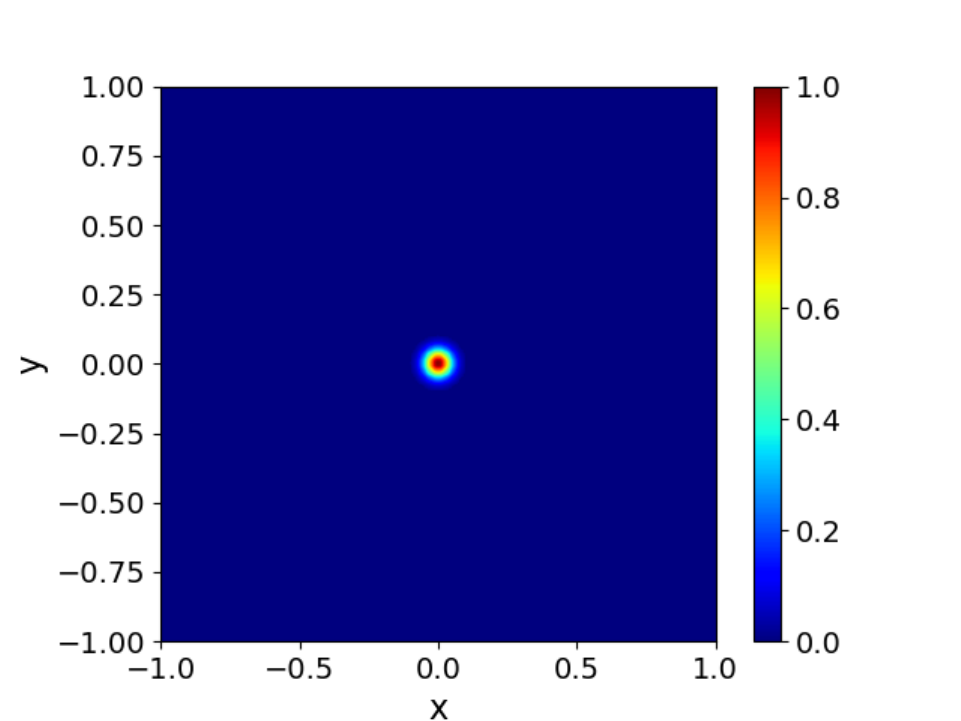}
	\end{minipage}
    }%
    \subfigure[t=1]{
	\begin{minipage}{0.49\linewidth}
		\centering
		\includegraphics[width=0.9\linewidth]{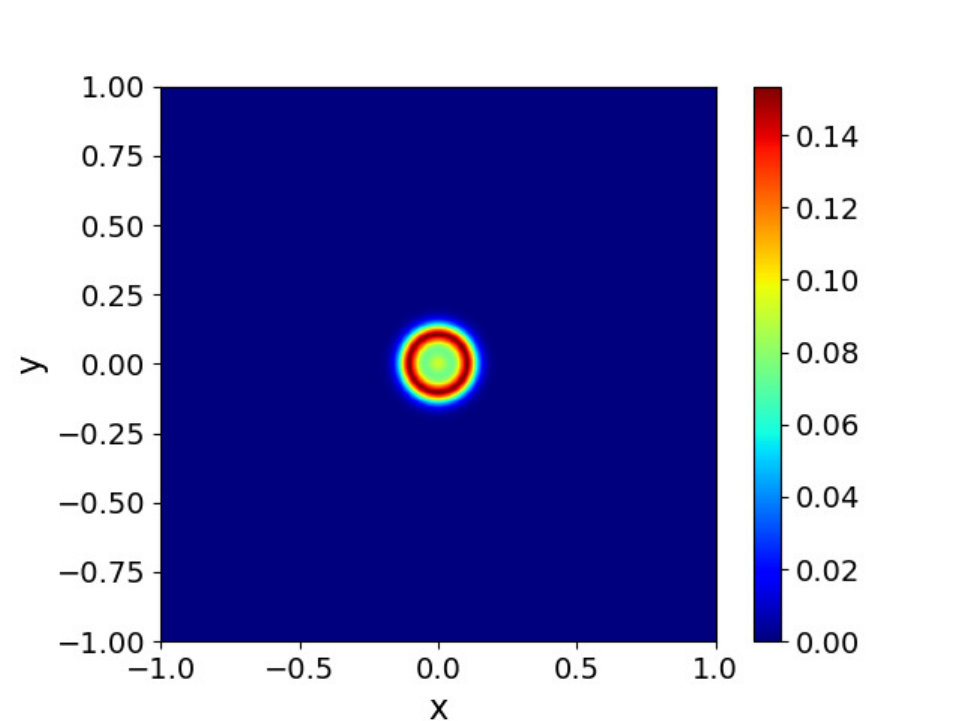}
	\end{minipage}
	}%
    \\ 
	\subfigure[t=2]{
	\begin{minipage}{0.49\linewidth}
		\centering
		\includegraphics[width=0.9\linewidth]{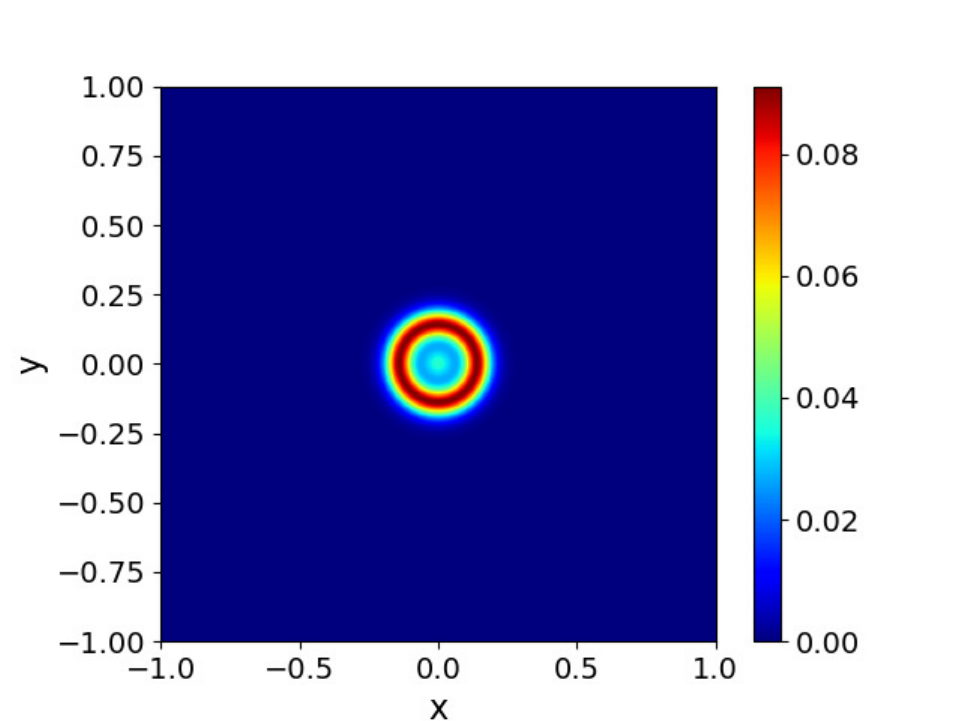}
	\end{minipage}
    }%
    \subfigure[t=4]{
	\begin{minipage}{0.49\linewidth}
		\centering
		\includegraphics[width=0.9\linewidth]{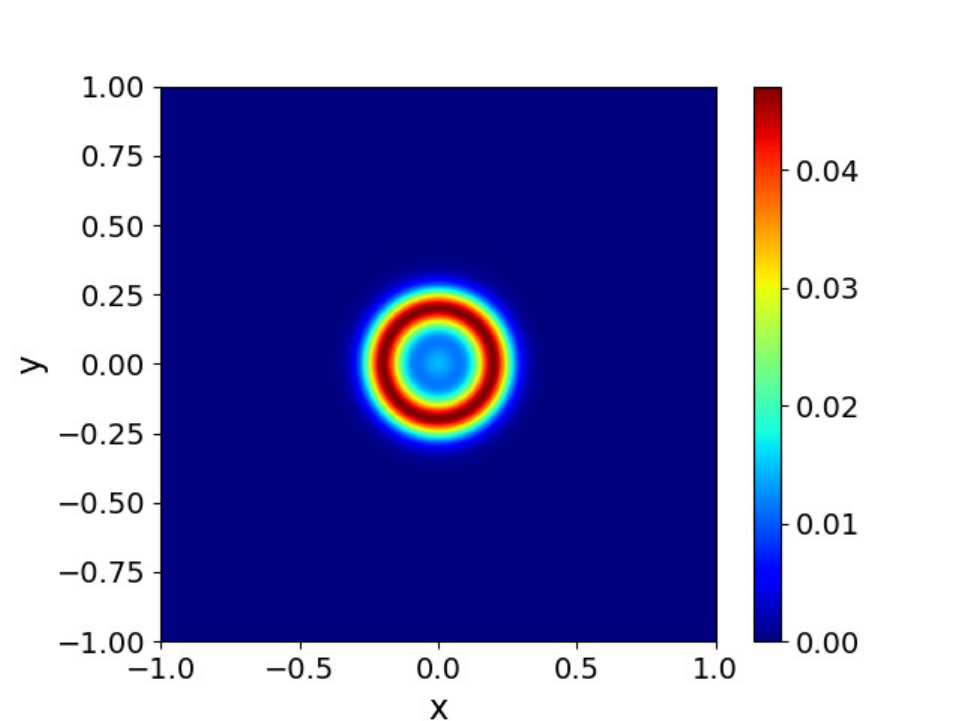}
	\end{minipage}
    }%
\caption{A numerical simulation of the system \eqref{cancer system}, with constant tumor cell diffusion, reveals the spatio-temporal dynamics of the tumor invasion process. The figure shows the emergence of a ring of cells that breaks away from the primary tumor mass and invades deeper into the ECM.}
\label{four snapshots}
\end{figure}

To demonstrate the spatio-temporal dynamics of the cancer cell invasion HAD model, we first repeat a 2D cancer cell spreading numerical experiment of  \cite{anderson2000mathematical}, which was derived from the FDM approximation of the system \eqref{cancer system}. The four snapshots shown in Figure \ref{four snapshots} illustrate the temporal progression of the tumor cell density distribution, with the first sub-figure representing the initial conditions:

Fig.\ref{four snapshots} shows four snapshots in time of the tumor cell density distribution, with the first figure representing the initial data:
\begin{equation}
\rho(x,y,0)=e^{-r^2/\epsilon},r^2=x^2+y^2,r\in[0,0.1]
\end{equation}
\begin{equation}
f(x,y,0)=1-0.5\, \rho(x,y,0)
\end{equation}
\begin{equation}
m(x,y,0)=0.5\,\rho(x,y,0),
\end{equation}
where we set the same parameters as \cite{anderson2000mathematical} in the system \eqref{cancer system}: 
\begin{equation}\label{parameter}
d_n = d_m = 0.001, \quad \gamma = 0.005, \quad \eta = 10, \quad \alpha = 0.1, \quad \beta = 0, \quad \epsilon = 0.0025.
\end{equation}

In the tumor cell equation \eqref{cancer system}, the absence of terms for cell birth and death, combined with the application of zero flux boundary conditions, ensures that the total number of cells remains constant. The conservation of cell numbers allows us to verify the precision of the FDM. To quantify the deviation from the expected conservation,  we define the error at time $t=T$:

\begin{equation}\label{conservation}
\text{Error$_{t=T}$} = \frac{\sum_i(\rho_{i,t=T}) - \sum_i (\rho_{i,t = 0})}{\sum_i (\rho_{i,t = 0})}
\end{equation}
where $\rho_i$ refers to the density of the cell within the $i$-th grid of FDM. It has demonstrated an accuracy within $0.01\%$, indicating high reliability in the numerical simulation. We observe that the main body of the tumor invades slowly. At the forefront, a high-density cell zone emerges, which subsequently detaches to form an independent circular cluster of cells that penetrates deeper into the ECM.

Our experiments here and below are all carried out on the HPC2021 system at the University of Hong Kong, equipped with 16-core Intel Xeon 6226R processors, and an NVIDIA Tesla V100 32GB SXM2 GPU.

\subsection{Comparing radial/FDM/SIPF methods in 3D}
The goal of this section is to generalize the model to a 3D spatial domain, enabling a more detailed exploration of the spatio-temporal evolution of the system.
We set the initial conditions:
\begin{equation}\label{3D initial condition}
\begin{aligned}
\rho(x,y,z,0)&=e^{-r^2/\epsilon},r^2=x^2+y^2+z^2,r\in[0,0.1], \\
f(x,y,z,0)&=1-0.5\, \rho(x,y,z,0), \\
m(x,y,z,0)&=0.5\, \rho(x,y,z,0).
\end{aligned}
\end{equation}

Unless otherwise stated, the parameter values utilized in the subsequent simulations were the same as those employed in the previous 2D experiments, as given by \eqref{parameter}.
Here we know that the condition is the radially symmetric case, we can write $\rho(x,y,z,t)=\rho(r,t)$, $f(x,y,z,t)=f(r,t)$, $m(x,y,z,t)=m(r,t)$ and simplify the equation as:
\begin{equation}
\rho_t = d_n\, \left (\frac{\partial^2 \rho}{\partial r^2}+\frac{2}{r}\frac{\partial \rho}{\partial r}\right ) - \gamma\, \left (\frac{\partial \rho}{\partial r}\frac{\partial f}{\partial r} + \rho\cdot(\frac{\partial^2 f}{\partial r^2}+\frac{2}{r}\frac{\partial f}{\partial r})\right ); 
\end{equation}

\begin{equation}
f_t = -\eta m f; 
\end{equation}

\begin{equation}
m_t = d_m\, \left (\frac{\partial^2 m}{\partial r^2}+\frac{2}{r}\frac{\partial m}{\partial r}\right ) - \beta m + \alpha \rho.
\end{equation}

\begin{figure}[htbp]
	\centering
    \vspace{-0.2cm}
    \subfigtopskip=2pt
    \subfigure[t=1]{
	\begin{minipage}{0.49\linewidth}
		\centering
		\includegraphics[width=0.9\linewidth]{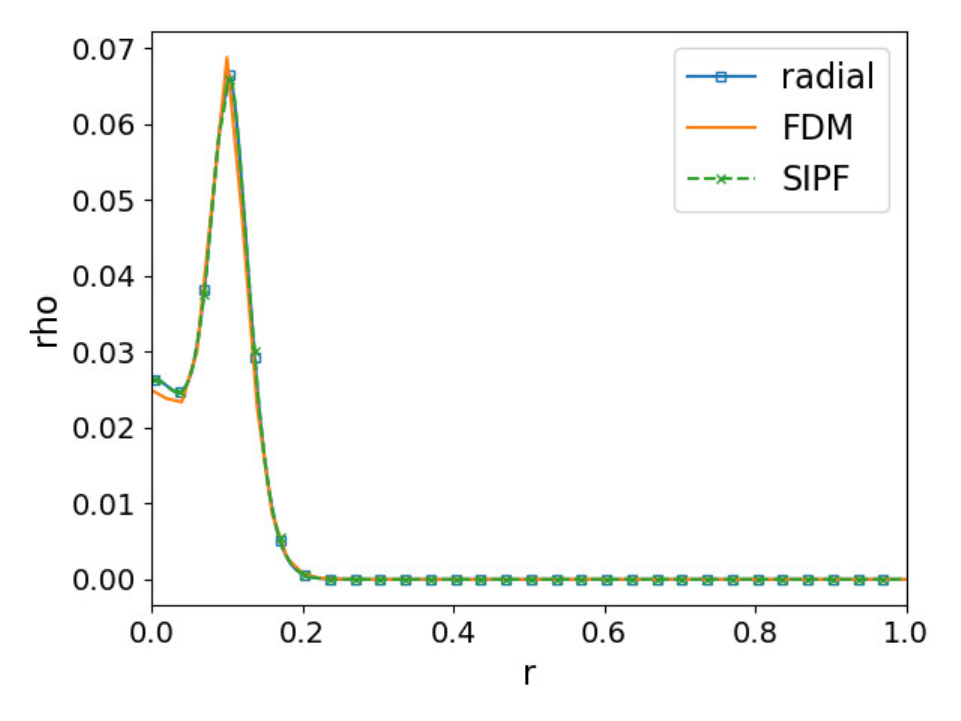}
	\end{minipage}
    }%
    \subfigure[t=2]{
	\begin{minipage}{0.49\linewidth}
		\centering
		\includegraphics[width=0.9\linewidth]{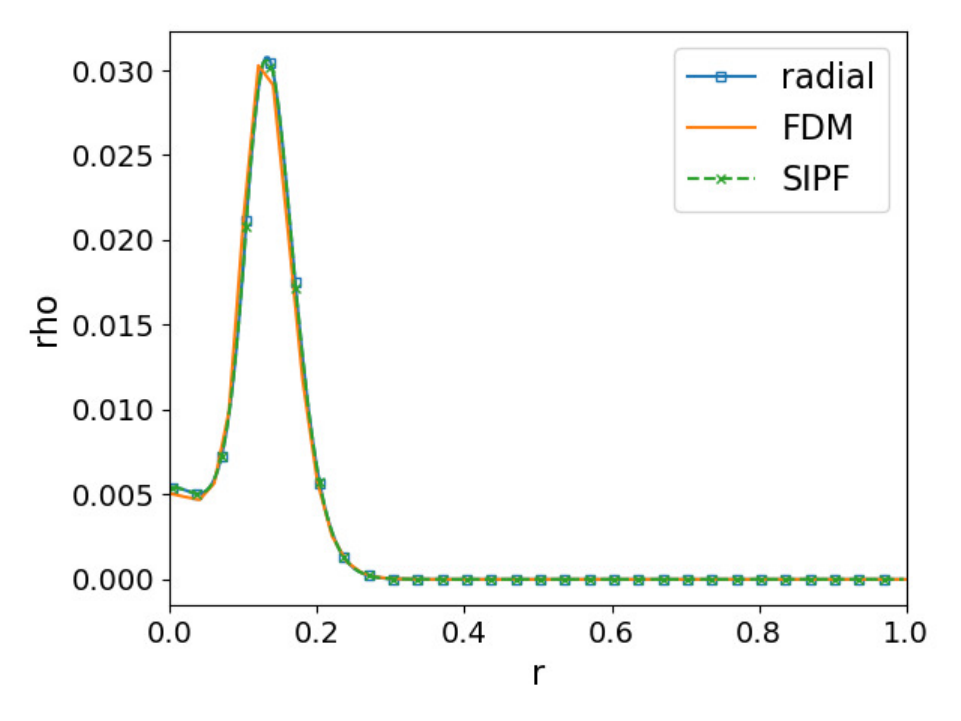}
	\end{minipage}
	}%
    \\ 
	\subfigure[t=3]{
	\begin{minipage}{0.49\linewidth}
		\centering
		\includegraphics[width=0.9\linewidth]{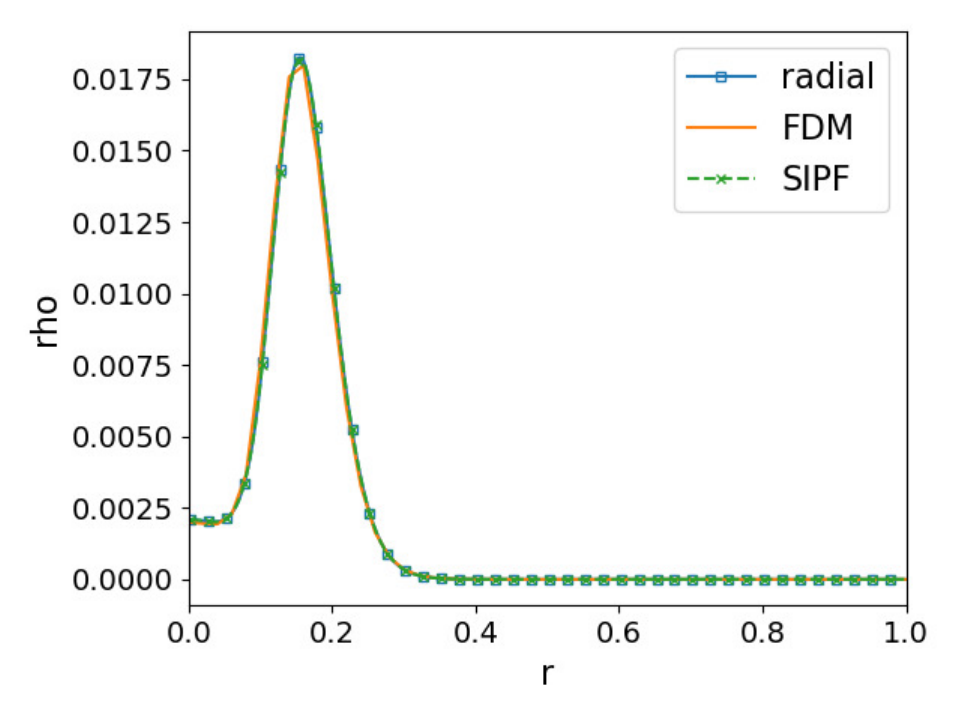}
	\end{minipage}
    }%
    \subfigure[t=4]{
	\begin{minipage}{0.49\linewidth}
		\centering
		\includegraphics[width=0.9\linewidth]{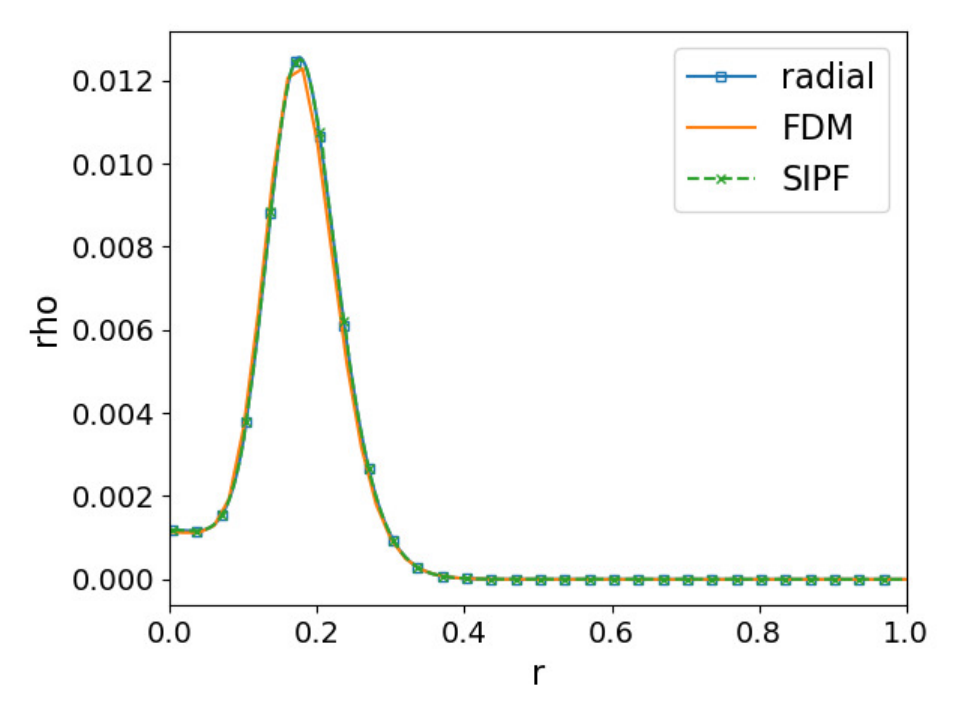}
	\end{minipage}
    }%
\caption{3D numerical solutions of the system \eqref{cancer system} with constant tumor cell diffusion showing the cell density computed by radial, FDM, and SIPF.}
\label{3D FDM/SIPF/radial}
\end{figure}

We use a very fine mesh to compute the radial solution, which will be the reference solution in our numerical experiment. We compare the FDM and SIPF with radial solutions in this experiment. We conduct FDM numerical simulations on a uniform mesh with $\delta x = \delta y = \delta z= 1/101$, with a time step \(\delta t = 10^{-2}\), radial 1D simulations on a uniform mesh with $\delta r = 1/301$, \(\delta t = 10^{-3}\). For the SIPF method, we discretize the MDE concentration \(m\) using $H=24$ Fourier basis in each spatial dimension and approximate the distribution \(\rho\) with $P=10,000$ particles. We simulate the evolution of \(m, f\) and \(\rho\) using Alg.\ref{the whole algorithm}, with a time step \(\delta t = 10^{-2}\). Fig.\ref{3D FDM/SIPF/radial} presents 1D slices of the results from the FDM and SIPF methods, depicting the temporal progression of tumor cell invasion into the host tissue. The SIPF method demonstrates higher accuracy than FDM, particularly at the peak values.

We also give three tables to further demonstrate two methods. To measure the convergence rate, we compute the solutions of FDM and SIPF on different conditions and compare the obtained relative $L^2$ error of $m$ with the reference solution, computed by the proposed radial case on a uniform mesh with $\delta r = 1/801$ and \(\delta t = 10^{-3}\). Taking SIPF as an example, we convert 3D spatial domain data to a 1D radial representation and compare it with the reference radial solution.
The spatial domain data, denoted as \( M \), is derived from the frequency domain data, which are the Fourier coefficients represented by \( \tilde{M}_\alpha \). Specifically, \( M \) is obtained by taking the real part of the result from the Inverse Fast Fourier Transform (IFFT) applied to \( \tilde{M}_\alpha \). We define \( n \) bins with edges \( r_i \) for \( i = 0 \) to \( n \).
The relative \( L^2 \) error between the SIPF mean and the reference radial solution is defined as:
\begin{equation}
\text{Relative $L^2$ Error} = \frac{\sqrt{\sum_{i} (M_i - R_i)^2}}{\sqrt{\sum_{i} (R_i)^2}}
\end{equation}
where \( M_i \) represents the SIPF mean in the \(i\)-th bin, and \( R_i \) denotes the radial reference value in the same bin.

Taking Table \ref{tab:FDM gridsize} as an example, the rate of convergence, denoted as \textit{Rate}, is defined through the formula:
\begin{equation}
\text{Rate} = \left| \frac{\log(\epsilon_{\text{prev}} / \epsilon_{\text{curr}})}{\log(\delta x_{\text{prev}} / \delta x_{\text{curr}})} \right|
\end{equation}
where:
\begin{itemize}
    \item \( \epsilon_{\text{prev}} \) is the relative \( L^2 \) error at the previous grid size,
    \item \( \epsilon_{\text{curr}} \) is the relative \( L^2 \) error at the current grid size,
    \item \( \delta x_{\text{prev}} \) and \( \delta x_{\text{curr}} \) are the spatial step sizes at the previous and current grid sizes, respectively.
\end{itemize}

Similarly, we can define the ratio in this way. According to Table \ref{tab:FDM gridsize}, the FDM method not only proves to be inaccurate but also time-consuming. As the grid size increases, the computational runtime of FDM escalates significantly. However, the accuracy of the numerical method does not improve commensurately with the finer grid, failing to justify the substantial increase in runtime. Table \ref{tab: SIPF dt} illustrates the variations in computational runtime and relative $L^2$ error with changes in the time step $\delta t$ for the SIPF method. Unlike FDM, SIPF is not as time-intensive, and larger $\delta t$ values still maintain commendable accuracy. Additionally, Table \ref{tab: SIPF particle} shows that increasing the number of particles significantly impacts the runtime. 

\begin{figure}[htbp]
    \centering
    \subfigtopskip=2pt
    \subfigbottomskip=2pt
    \subfigure[vs. time step $\delta t$ on log-scale]{
        \begin{minipage}{0.475\linewidth}
            \centering
            \includegraphics[width=\linewidth]{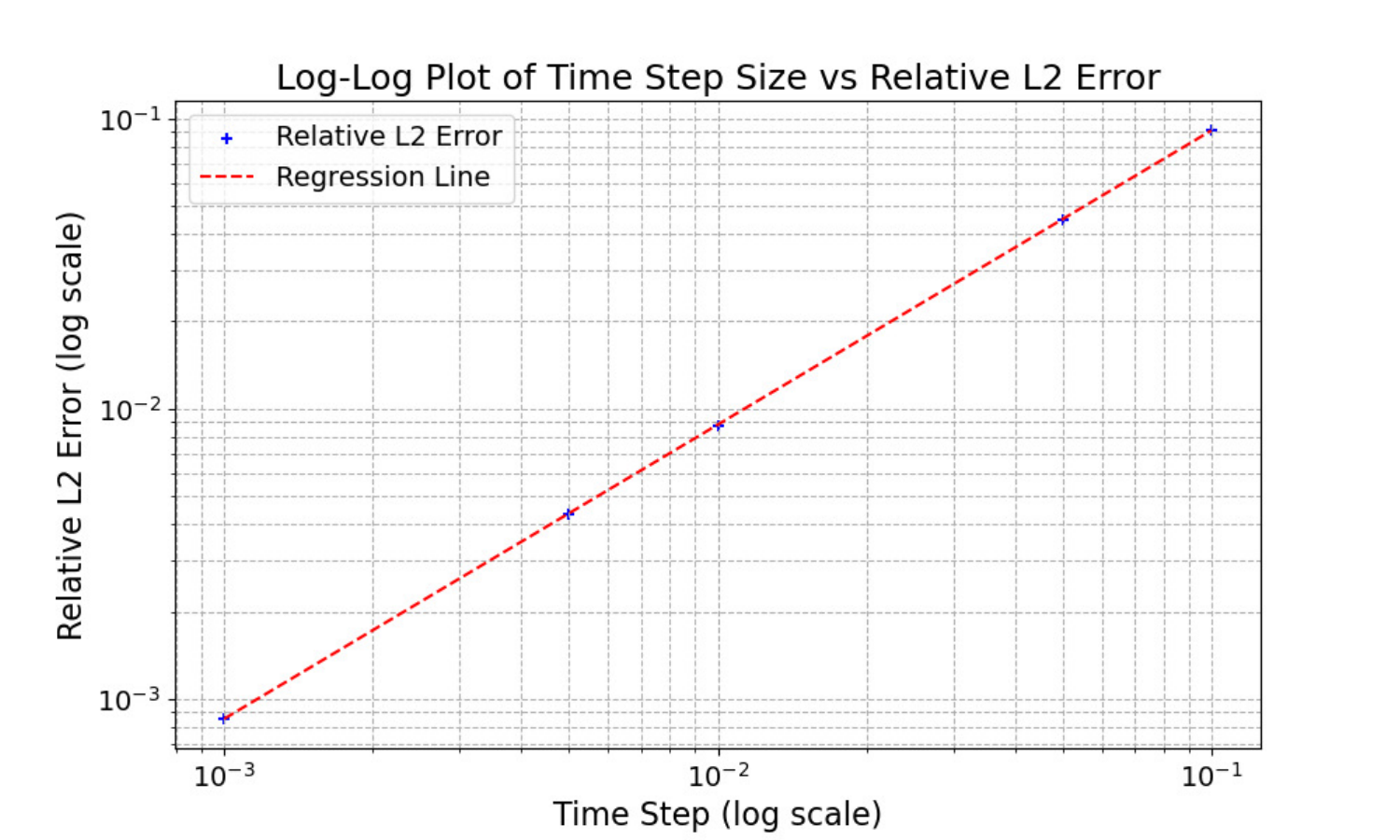}
            \label{fig:SIPF_dt_L2error}
        \end{minipage}
    }
    \subfigure[vs. particle number $P$ on log-scale]{
        \begin{minipage}{0.475\linewidth}
            \centering
            \includegraphics[width=\linewidth]{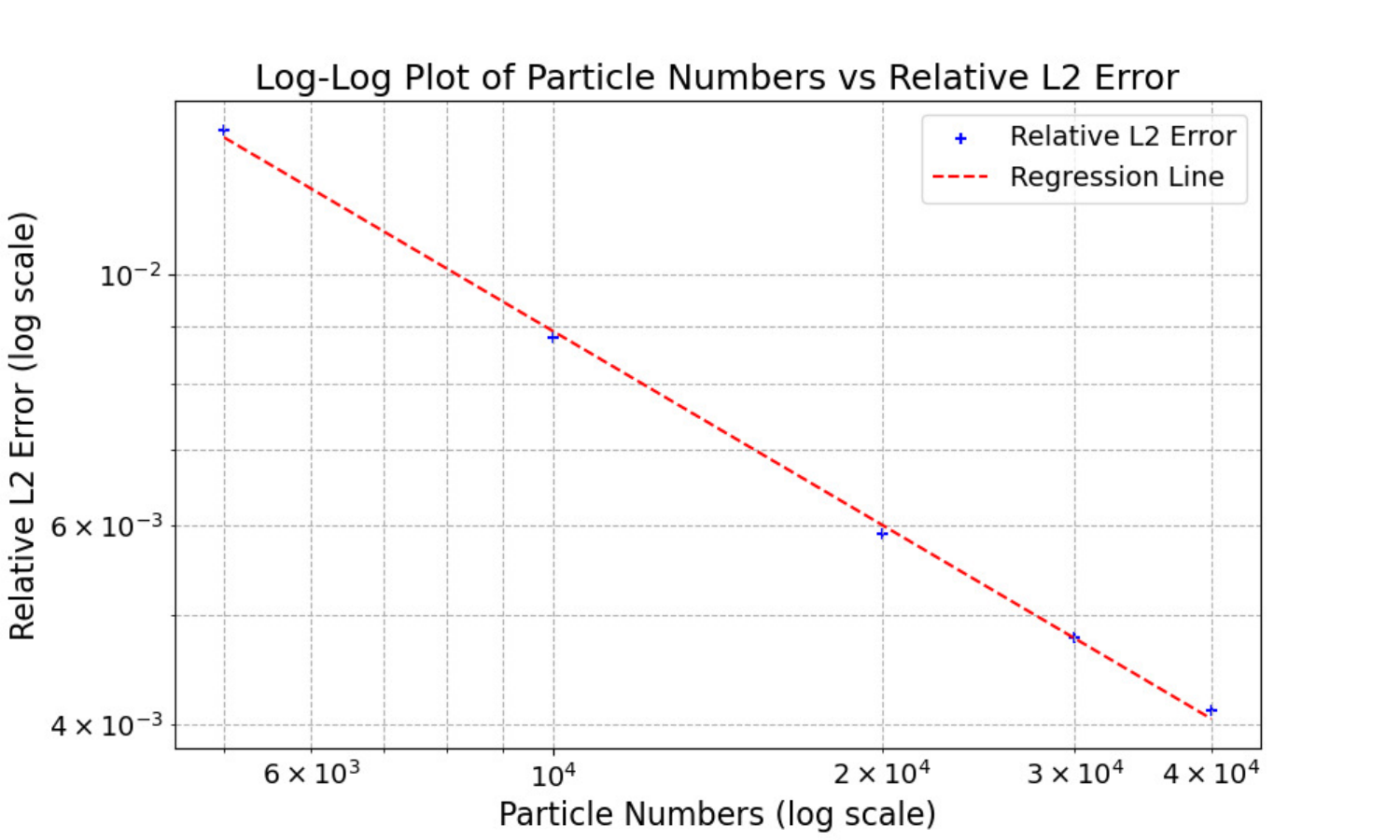}
            \label{fig:SIPF_P_L2error}
        \end{minipage}
    }
    \vspace{2em}
    \subfigure[vs. Fourier mode $H$ on semi-log-scale]{
        \begin{minipage}{0.475\linewidth}
            \centering
            \includegraphics[width=\linewidth]{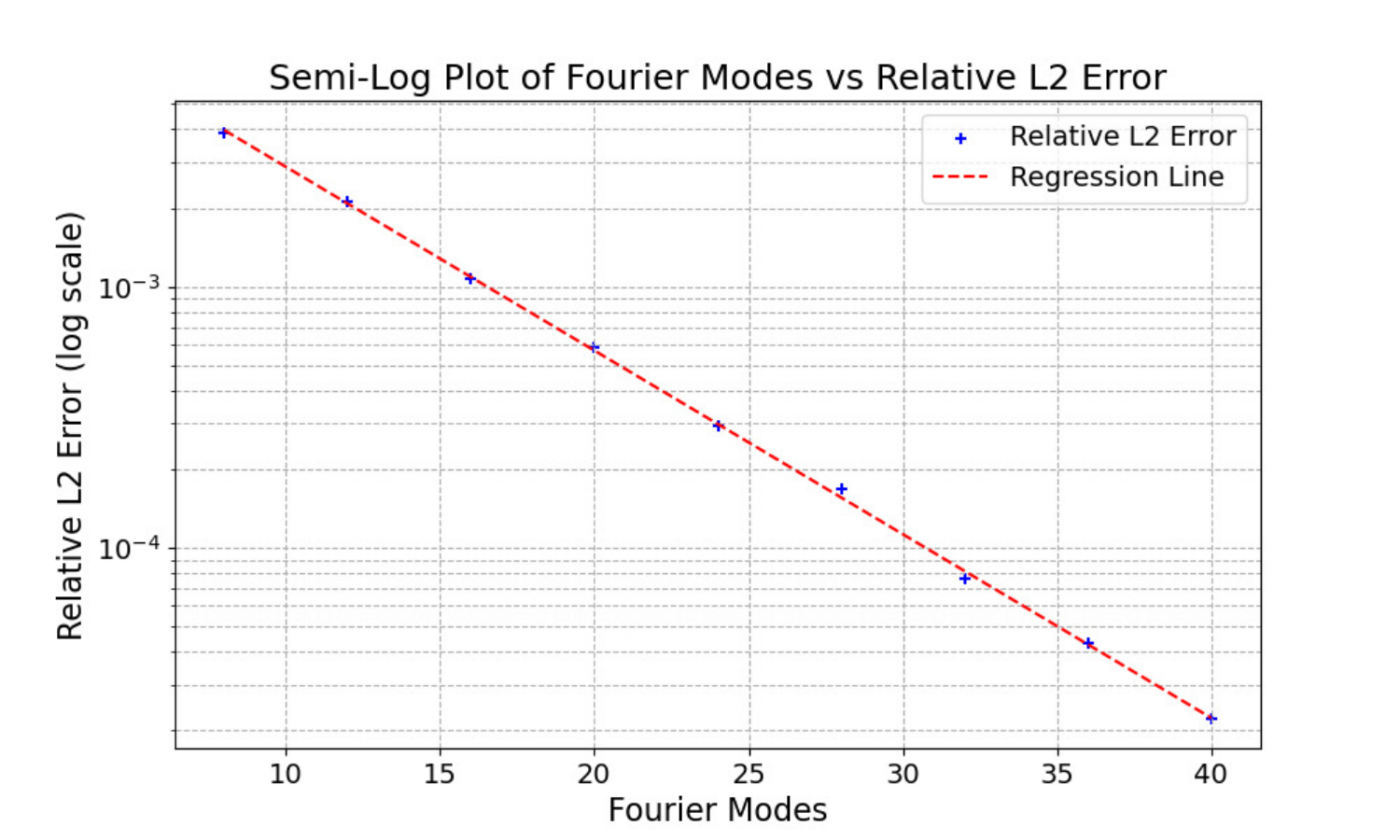}
            \label{fig:SIPF_H_L2error}
        \end{minipage}
    }
    \vspace{-1.5em}
    \caption{3D relative $L^2$ errors of $m$ in SIPF (radial solution being the reference).}
    \label{SIPF_L2error}
\end{figure}

In Fig.\ref{SIPF_L2error}, we compute the relative $L^2$ error of the MDE concentration $m$ at the final time $T=4$ for different time step $\delta t$, particle number $P$, and Fourier mode $H$. In addition, we fit the slope of the error versus $\delta t$ in the logarithmic scale and find $e(\delta t) = \mathcal{O}(\delta t^{1.0130})$, indicating the algorithm being approximately first-order in time. Furthermore, by fitting the slope of the error versus $P$ in the logarithmic scale, we find that $e(P) = \mathcal{O}(P^{-0.5587})$. To provide a clearer picture of the convergence of the MDE concentration $m$ versus Fourier mode $H$, we plot in Fig.\ref{fig:SIPF_H_L2error} the errors in the semi-log scale. This plot indicates an exponential convergence rate $\mathcal{O}(e^{-0.1608H})$. In terms of accuracy, measured by the relative $L^2$ error, there is a clear improvement as the number of Fourier modes increases. Experiments indicate that when we set particle number $P$ to be 10,000, time step $\delta t$ to be 0.01, and Fourier mode $H$ to be 24, there is a good trade-off between accuracy and computational time. The following SIPF algorithm adopts this configuration with no specific mention. 

\begin{table}[!htbp]
\centering
\begin{tabular}{ccccc}  
\toprule
FDM Grid & Run time(s) & Ratio & Relative $L^2$ Error & Rate\\
\midrule 
$21 \times 21 \times 21$ & 15.54 &  & 1.1430 &  \\
$41 \times 41 \times 41$ & 132.87 & 3.09 & 0.2808 & 2.02\\
$61 \times 61 \times 61$ & 465.08 & 3.09 & 0.1253 & 1.99\\
$81 \times 81 \times 81$ & 1126.01 & 3.07 & 0.0694 & 2.05\\
$101 \times 101 \times 101$ & 2238.85 & 3.08 & 0.0447 & 2.01\\
\bottomrule
\end{tabular}
\caption{3D run time and relative $L^2$ error of FDM vs. grid size (at $\delta t=0.01$).}\label{tab:FDM gridsize}
\end{table}

\begin{table}[!htbp]
\centering
\begin{tabular}{ccccc}  
\toprule
dt(s) & Run time(s) & Ratio & Relative $L^2$ Error & Rate\\
\midrule
0.1 & 10.37 &  & 9.12E-02 & \\ 
0.05 & 18.06 & 0.80 & 4.50E-02 & 1.02\\ 
0.01 & 66.98 & 0.81 & 8.78E-03 & 1.02\\
0.005 & 117.42 & 0.81 & 4.36E-03 & 1.01\\ 
0.001 & 413.01 & 0.79 & 8.59E-04 & 1.01\\
\bottomrule
\end{tabular}
\caption{3D run time and relative $L^2$ error of SIPF vs. $\delta t$ (at $P=10000$).}\label{tab: SIPF dt}
\end{table}

\begin{table}[!htbp]
\centering
\begin{tabular}{ccccc}  
\toprule
Particle Numbers & Run time(s) & Ratio & Relative $L^2$ error & Rate\\
\midrule
5000 & 29.70 &  & 1.34E-02 & \\
10000 & 66.98 & 0.95 & 8.78E-03 & 0.61\\
20000 & 358.44 & 2.42 & 5.89E-03 & 0.58\\
30000 & 1129.15 & 2.83 & 4.77E-03 & 0.52\\
40000 & 1995.86 & 1.98 & 4.11E-03 & 0.51\\
\bottomrule
\end{tabular}
\caption{3D run time and relative $L^2$ error of SIPF vs. $P$ (at $\delta t=0.01$).}\label{tab: SIPF particle}
\end{table}

\subsection{Regime of Small Diffusion Coefficient}
\begin{figure}[htbp]
	\centering
    \vspace{-0.2cm}
    \subfigtopskip=2pt
    \subfigure[t=1]{
	\begin{minipage}{0.49\linewidth}
		\centering
		\includegraphics[width=0.9\linewidth]{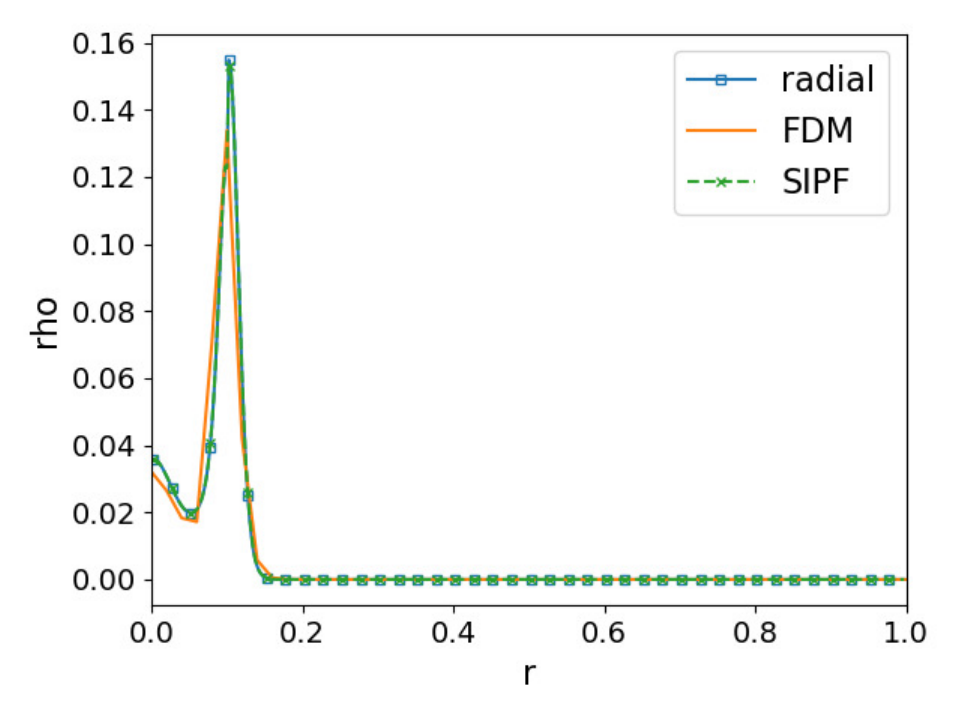}
	\end{minipage}
    }%
    \subfigure[t=2]{
	\begin{minipage}{0.49\linewidth}
		\centering
		\includegraphics[width=0.9\linewidth]{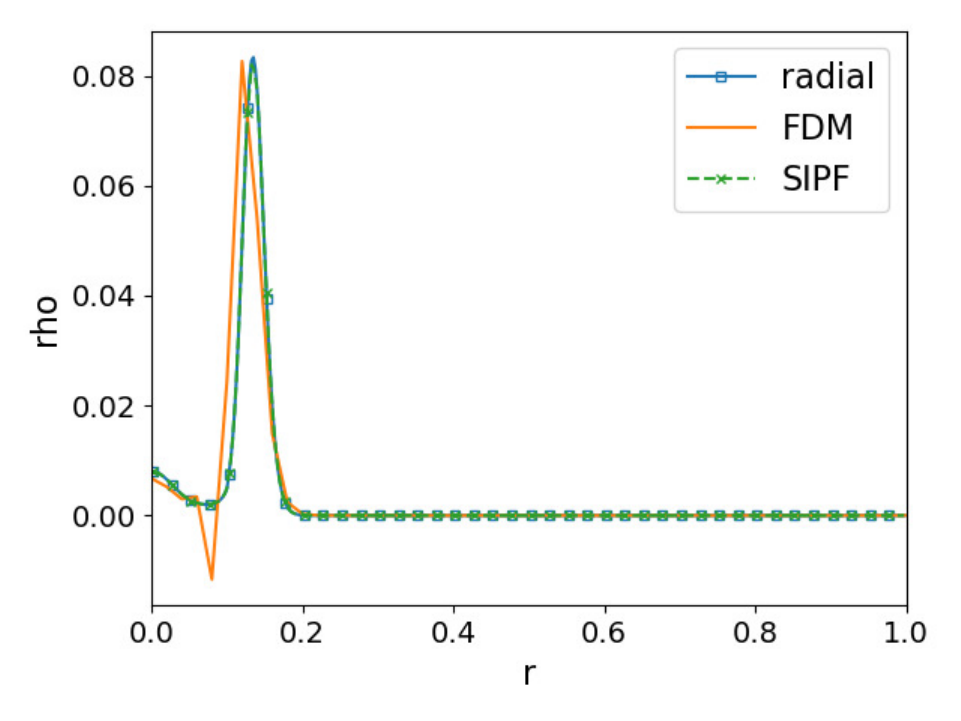}
	\end{minipage}
	}%
    \\ 
	\subfigure[t=3]{
	\begin{minipage}{0.49\linewidth}
		\centering
		\includegraphics[width=0.9\linewidth]{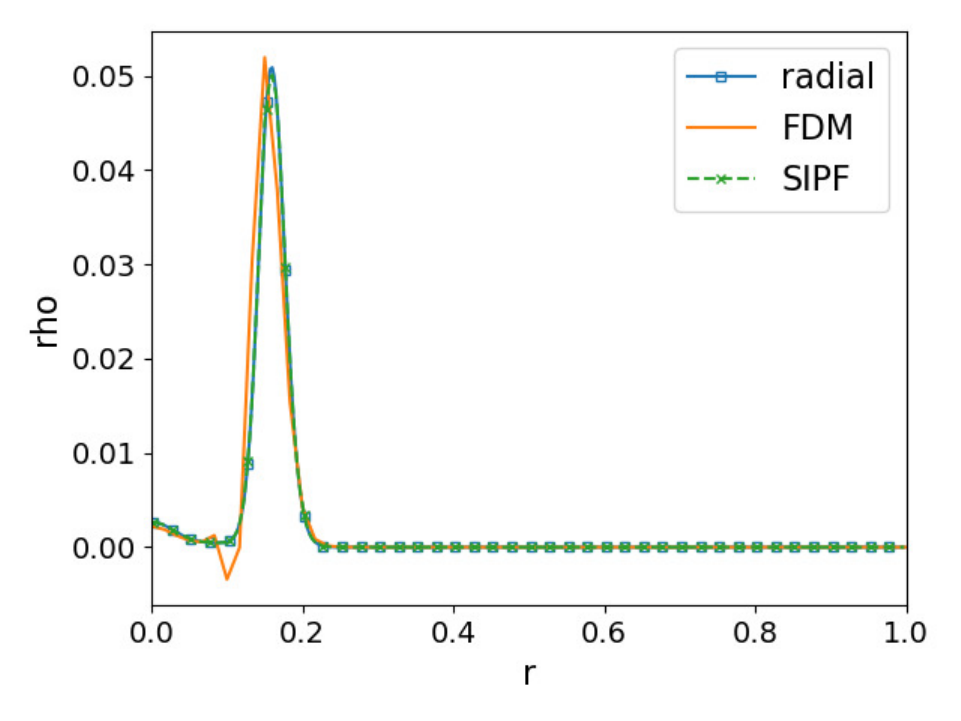}
	\end{minipage}
    }%
    \subfigure[t=4]{
	\begin{minipage}{0.49\linewidth}
		\centering
		\includegraphics[width=0.9\linewidth]{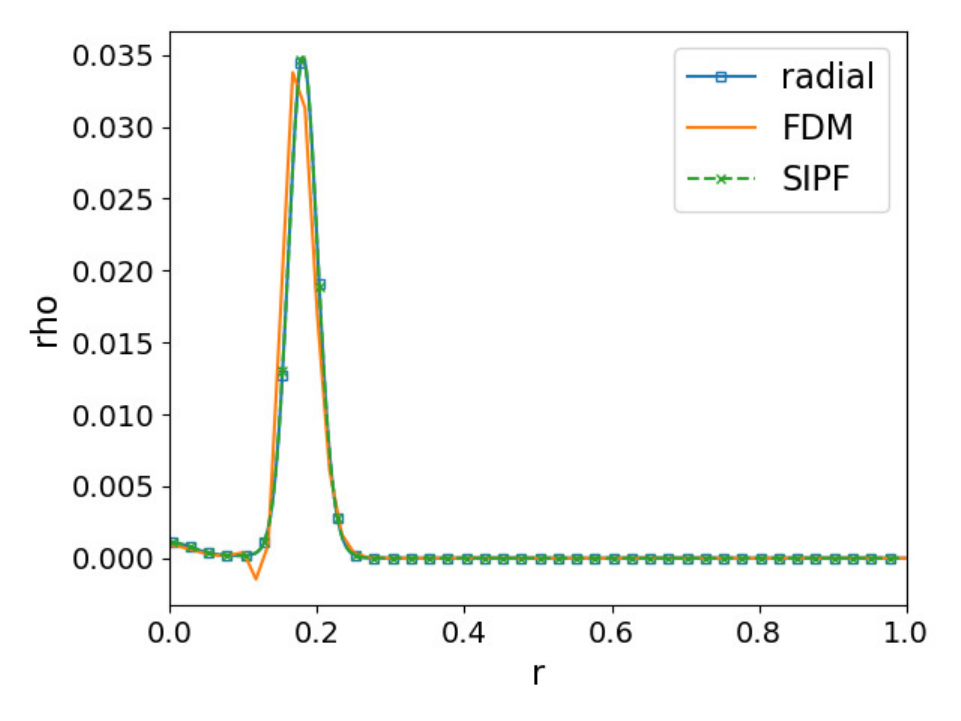}
	\end{minipage}
    }%
\caption{Comparing radial solutions, FDM (red) and SIPF (green) at $d_n=0.0002$, shows under-shoot (violating positivity) and inaccurate peak locations in FDM.}
\label{dn=0.0002}
\end{figure}

We change the diffusion coefficient $d_n$ of section 3.1, with other conditions and parameters unchanged. When $d_n$ becomes smaller, the FDM would be more expensive. In Fig.\ref{dn=0.0002}, we set $d_n=0.0002$. Compared with Fig.\ref{3D FDM/SIPF/radial}, we can see the simulation of Fig.\ref{dn=0.0002} is not good, as the diffusion coefficient $d_n$ decreases, the peak of tumor density becomes steeper. The FDM exhibits instability because it necessitates a very tight discretization to adequately resolve the peak, which in turn incurs a substantial time cost. However, the SIPF method remains stable and maintains high accuracy under these conditions.

\subsection{Two Clusters of Cancer Cells Evolution}
\begin{figure}[htbp]
	\centering
    \vspace{-0.2cm}
    \subfigtopskip=2pt
    \subfigure[t=0]{
	\begin{minipage}{0.49\linewidth}
		\centering
		\includegraphics[width=0.9\linewidth]{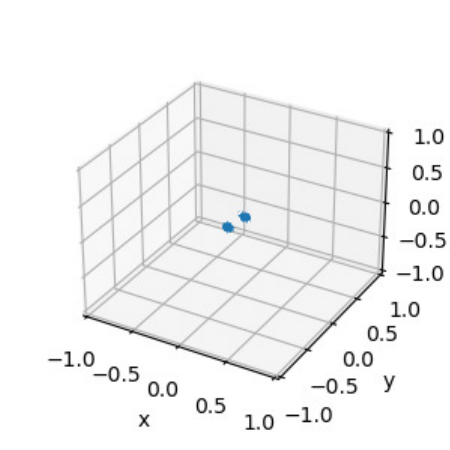}
	\end{minipage}
    }%
    \subfigure[t=0.3]{
	\begin{minipage}{0.49\linewidth}
		\centering
		\includegraphics[width=0.9\linewidth]{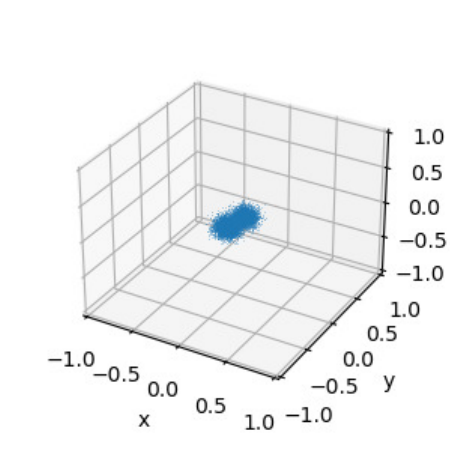}
	\end{minipage}
	}%
    \\ 
	\subfigure[t=0.7]{
	\begin{minipage}{0.49\linewidth}
		\centering
		\includegraphics[width=0.9\linewidth]{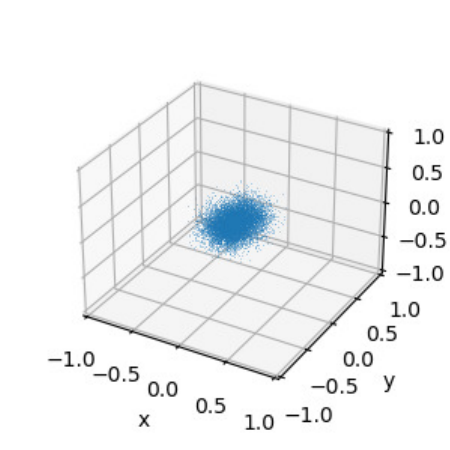}
	\end{minipage}
    }%
    \subfigure[t=1]{
	\begin{minipage}{0.49\linewidth}
		\centering
		\includegraphics[width=0.9\linewidth]{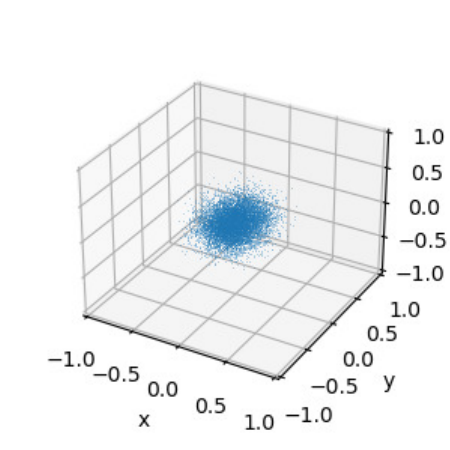}
	\end{minipage}
    }%
\caption{Two clusters of cancer cells merge into a larger single cluster and spread further.}
\label{diffusion fusion}
\end{figure}

The aforementioned experiments all satisfy the symmetry of the initial data, in this subsection, we investigate the behaviors from non-radial initial data. We demonstrate that the SIPF algorithm is equally applicable to asymmetric situations.
We show a SIPF simulation on two clusters of cancer cells spreading dynamics in 3D. The initial condition is:
\begin{equation}
\rho_1(x,y,z,0)=e^{-r_1^2/\epsilon},r_1^2=(x-a)^2+(y-b)^2+(z-c)^2,r\in[0,0.1]
\end{equation}
\begin{equation}
\rho_2(x,y,z,0)=e^{-r_2^2/\epsilon},r_2^2=(x-d)^2+(y-e)^2+(z-f)^2,r\in[0,0.1]
\end{equation}
\begin{equation}
\rho(x,y,z,0)=\rho_1(x,y,z,0)+\rho_2(x,y,z,0).
\end{equation}
Here we choose $a=b=c=0.1,d=e=f=-0.1$ and display the fusion/spreading process of these two clusters in Fig.\ref{diffusion fusion}. Two clusters of cells diffuse outward over time, intersect, and subsequently merge, continuing their outward invasion. In the absence of any specified interactions between the two clusters of cells, the diffusion process remains analogous to that of a single cluster of cells.

\subsection{Comparing FDM and SIPF based on Integral Identities}
Given the initial condition \eqref{3D initial condition} in section 3.1, we have: 
\[\int_{\Omega} m(X, 0) \, dX = 2\pi\int_{0}^{0.1}e^{-\frac{r^2}{\epsilon}} r^2\,dr = \frac{1}{2}\int_{\Omega} \rho(X, 0) \, dX, \; \beta=0,\]
according to \eqref{integral m}, and the following relationship holds:
\begin{equation}\label{integration m}
\int_{\Omega} m(X, t) = 4\pi\alpha t \int_{0}^{0.1}e^{-\frac{r^2}{\epsilon}} r^2\,dr + 2\pi\int_{0}^{0.1}e^{-\frac{r^2}{\epsilon}}r^2\,dr.
\end{equation} 

We compute $\int_{\Omega} m(X, t)$ using the radial solution, FDM, and SIPF methods, and compare the results with the reference value in \eqref{integration m}.
Given Eq.\eqref{lnf}, Eq.\eqref{integration m} and \[\int_{\Omega} \ln f_0\, dX= 4\pi \int_0^{0.1} \ln\left(1 - 0.5 e^{-r^2/\epsilon}\right) r^2 \, dr,\]
we have that $f$ satisfies
\begin{align}\label{eq:integration_lnf}
\int_{\Omega}  \ln f(X,t) &= 4\pi \int_0^{0.1} \ln\left(1 - 0.5 e^{-r^2/\epsilon}\right) r^2 \, dr \nonumber \\
&-\eta \left (2\pi\alpha t^2 \int_{0}^{0.1}e^{-\frac{r^2}{\epsilon}} r^2\,dr+ 2\pi t\int_{0}^{0.1}e^{-\frac{r^2}{\epsilon}}r^2\,dr\right ).  \end{align}

\begin{figure}[htbp]
	\centering
    \vspace{-0.2cm}
    \subfigtopskip=2pt
    \subfigure[$m$]{
	\begin{minipage}{0.49\linewidth}
		\centering
		\includegraphics[width=0.9\linewidth]{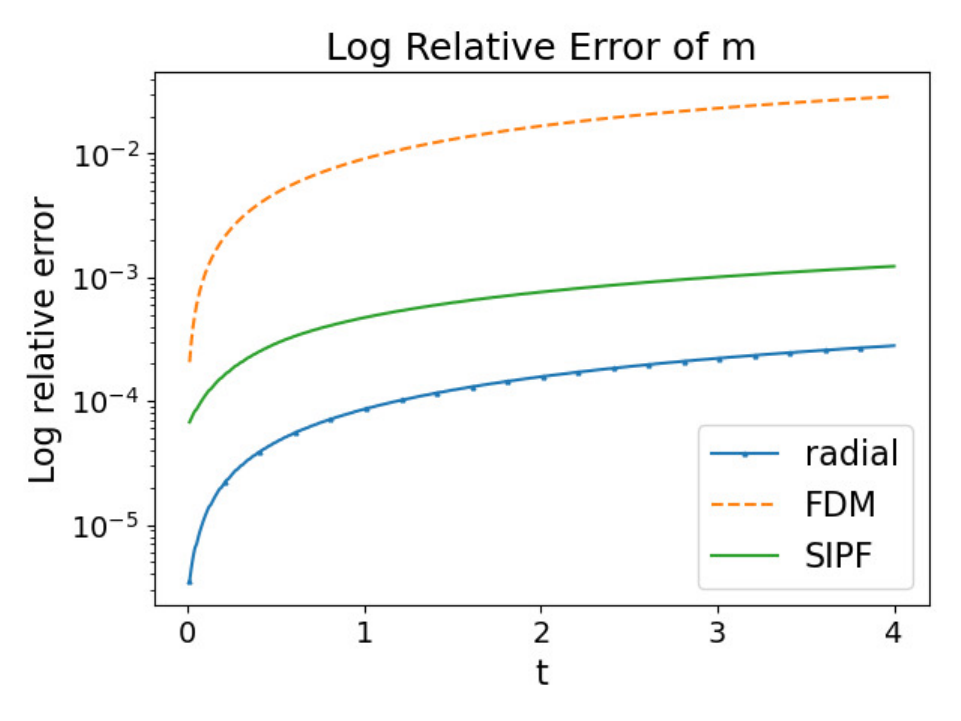}
	\end{minipage}
    }%
    \subfigure[$\ln \, f$]{
	\begin{minipage}{0.49\linewidth}
		\centering
	 \includegraphics[width=0.9\linewidth]{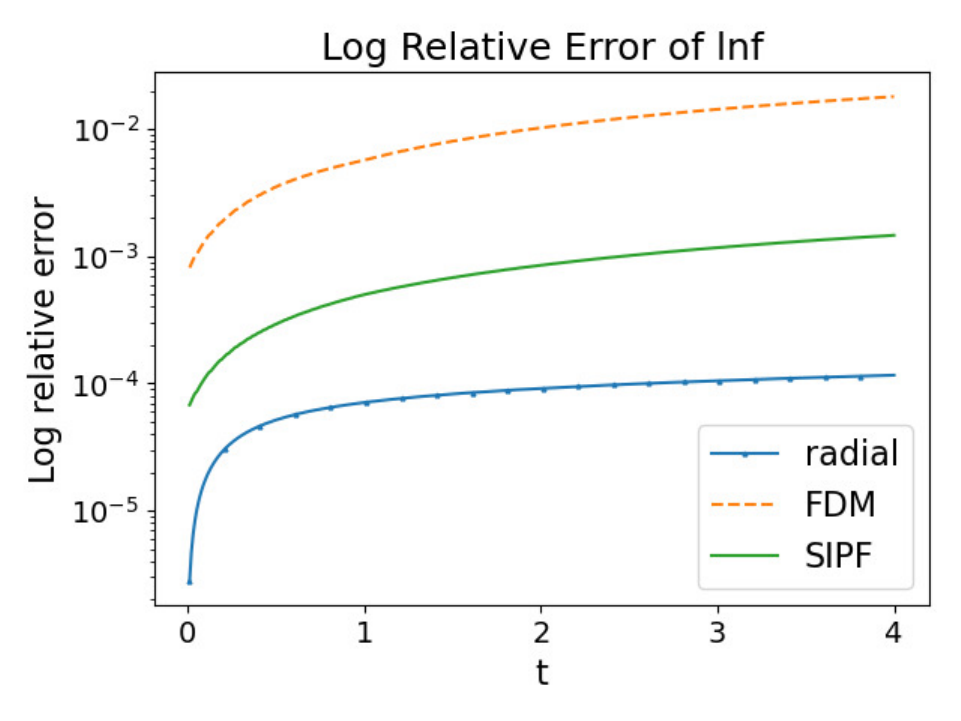}
	\end{minipage}
	}%
\caption{Relative error comparison in computing the integrals of Eqs.\eqref{integration m}\eqref{eq:integration_lnf} by radial solution (blue), FDM (red), and SIPF (green).}
\label{Relative Error}
\end{figure}

We define the relative \(L^2\) error of \(m\) at time \(T\), which involves several calculations for different methods of approximation or measurement.  First, we calculate the reference value of \(m\), denoted as \(m_r\), which is the integral of \(m(X, T)\) over the domain \(\Omega\) in Eq.\eqref{integration m}. For the SIPF method, \(m\) at time \(T\) is computed as \(m^\text{SIPF} = \alpha_{T;0,0,0} \cdot L^3\), where \(\alpha\) is a Fourier coefficient defined previously (see Eq.\eqref{m Fourier}) and \(L\) is a characteristic length scale. In the radial method, the radial coordinate $r$ is discretized into bins with widths $\delta r$ (the radial step size), and values $m$ corresponding to these discrete radii, the integral of \(m\) at time \(T\) is approximated as a sum: \(m^\text{radial} = \sum_i (m_{i,t=T}) \cdot 4\pi r_{i}^2 \cdot \delta r\), where $m_i$ is the value of $m$ at the $i$-th radial position, $r_i$ is the radius at the $i$-th position. For the FDM, we compute \(m\) at time \(T\) as \(m^\text{FDM} = \sum_j (m_{j,t=T}) \cdot (\delta x)^3\), with \(\delta x\) being the spatial step size used in the FDM. With these calculations, the relative error of \(m\) at time \(T\), denoted as \(Error_m\), is defined mathematically as:
\begin{equation}\label{m relative error}
\text{Error}_{\text{m}} = \frac{|m_n - m_r|}{|m_r|}
\end{equation}
where \(m_n\) represents any of the numerically computed values from the SIPF, radial, or FDM. Similarly, we can define the relative error of $\ln\, f$. Fig.\ref{Relative Error} shows the relative error of $m$ and $ln\, f$ in time. SIPF consistently demonstrates better performance compared to FDM. The relative error of both $m$ and $ln\, f$ is approximately an order of magnitude lower than FDM at the final time. Overall, SIPF outperforms FDM.


\section{Conclusion and Future Work}
In this paper, we developed the SIPF algorithm and demonstrated its efficacy and accuracy in computing the cancer cell invasion within the HAD system. The algorithm is recursive with no history dependence, and the field variable is computed by FFT. Due to the field variable (concentration) being smoother than the density, the FFT approach works with only dozens of Fourier modes. The spreading behavior of the density variable is resolved by 10k particles. We found that the regular implicit FDM is both time-consuming and inaccurate in 3D computation of tumor invasion. 

In future works, we will carry out a deep particle study \cite{DP_22,wang2024deepparticle} based on the data generated from SIPF simulations here and explore more complex models of tumor invasion to better capture the growth dynamics. We will incorporate the oxygen supply into the existing system to enable more precise computations \cite{ANDERSON1998857}.
Furthermore, the coupled two-species cancer invasion haptotaxis model has practical significance in the real-world application and understanding of realistic tumor progression \cite{dai2023boundedness}, which we have only briefly discussed as a non-radial 3D case study and could be further explored.


\section*{Acknowledgements}
\noindent
ZW was partially supported by NTU SUG-023162-00001, and JX by NSF grant DMS-2309520. ZZ was supported by the National Natural Science Foundation of China  (Project 12171406), the Hong Kong RGC grant (Projects 17307921 and 17304324), the Outstanding Young Researcher Award of HKU (2020-21), Seed Funding for Strategic Interdisciplinary Research Scheme 2021/22 (HKU), and Seed Funding from the HKU-TCL Joint Research Centre for Artificial Intelligence.

\printbibliography

\end{document}